\documentclass[11pt,letterpaper,reqno]{amsart}

\usepackage{amsmath,amssymb,amsthm,enumitem, mathtools}

\usepackage{tikz}
\usetikzlibrary{arrows.meta}
\numberwithin{equation}{section}

\usepackage[colorlinks=true,linkcolor=blue,citecolor=blue,urlcolor=blue]{hyperref}

\newcommand{\R}{\mathbb R}
\newcommand{\Sph}{\mathbb S}
\newcommand{\Cone}{\mathcal C}
\newcommand{\M}{\mathcal M}
\newcommand{\E}{\mathcal E}
\newcommand{\dist}{\operatorname{dist}}
\newcommand{\sgn}{\operatorname{sgn}}
\newtheorem{theorem}{Theorem}[section]
\newtheorem{lemma}[theorem]{Lemma}
\newtheorem{proposition}[theorem]{Proposition}
\newtheorem{corollary}[theorem]{Corollary}
\theoremstyle{definition}
\newtheorem{definition}[theorem]{Definition}

\title[Uniqueness of Allen--Cahn global minimizers]
{A uniqueness theorem for Allen--Cahn global minimizers
asymptotic to the Simons cone}

\author[Sanghoon Lee]{Sanghoon Lee}
\address[Sanghoon Lee]{Research Institute of Mathematics,
Seoul National University, 1 Gwanak-ro, Gwanak-gu,
Seoul 08826, Republic of Korea}
\email{lsh514@snu.ac.kr}

\author[Lorenzo Sarnataro]{Lorenzo Sarnataro}
\address[Lorenzo Sarnataro]{University of Toronto, Department of Mathematics, 40 St George Street, Toronto, ON M5S 2E4, Canada}
\email{lorenzo.sarnataro@utoronto.ca}

\date{\today}

\begin{document}

\begin{abstract}
For every integer $m\ge4$, we prove that the saddle-shaped solution
of the Allen--Cahn equation in $\R^{2m}$ constructed by Cabr\'e and
Terra \cite{CabreTerra2009} is a global minimizer. Moreover, up to sign, 
it is the unique global minimizer whose nodal set at infinity is a normal 
graph over the Simons cone with a height that decays faster than 
that of the Hardt--Simon leaves.
\end{abstract}

\maketitle

\section{Introduction}\label{sec:intro}

\noindent In this paper we shall consider the minimizing properties of certain bounded entire solutions of the Allen--Cahn equation
\begin{equation}\label{eq:AC}
	-\Delta u=u-u^3\qquad\text{in }\R^n.
\end{equation}
Solutions of \eqref{eq:AC} are critical points of the Allen--Cahn energy functional
\begin{equation}\label{eq:AC-energy}
	\E_\varepsilon(u;\Omega)\coloneqq\int_\Omega
	\left(\frac{\varepsilon}{2}|\nabla u|^2+\frac{W(u)}{\varepsilon}\right)d x
\end{equation}
with $\varepsilon=1$, potential $W:\R\to\R$ given by $W(s)=\frac14(1-s^2)^2$, and $\Omega\subset\R^n$ a bounded domain. 

We call $u$ a \emph{global minimizer} if $\E_1(u;\Omega)\leq\E_1(v;\Omega)$ for every smooth bounded domain $\Omega\subset\R^n$ and every $v\in C^1(\overline\Omega)$ with $v=u$ on $\partial\Omega$. We write $\E=\E_1$ for the unscaled energy. 

First proposed as a model for phase transitions, the Allen--Cahn energy $\E_\varepsilon$ has attracted significant mathematical interest due to its relation with the area functional and the theory of minimal hypersurfaces. Indeed, by~\cite{ModicaMortola1977,Sternberg88}, as $\varepsilon\downarrow0$, the functionals $\E_\varepsilon$ are known to $\Gamma$-converge to perimeter, and, by~\cite{Modica1979}, every sequence $u_\varepsilon:\Omega\to[-1,1]$ minimizing $\E_\varepsilon$ against compactly supported perturbations, with uniformly bounded energy, converges subsequentially in $L^1_{\mathrm{loc}}(\Omega)$ to a function of the form $\mathbf{1}_E-\mathbf{1}_{\Omega\setminus E}$, where $E$ locally minimizes perimeter in $\Omega$. 

By the interior regularity theory for area-minimizing boundaries~\cite{DeGiorgi1965, Almgren1966,Simons1968}, the limiting interface $\partial E$ is a smooth embedded hypersurface when $n\leq7$. Since area-minimizing hypersurfaces are locally modelled on area-minimizing cones, this regularity result can be seen as a consequence of the fact that every area-minimizing cone in $\R^n$ is flat for $n\leq 7$.

On the other hand, both regularity and rigidity fail in higher dimensions. Indeed, Bombieri, De Giorgi, and Giusti~\cite{BDG1969} proved that the 7-dimensional Simons cone
\[
\mathcal C^7=\{(x',x'')\in\R^4\times\R^4:|x'|=|x''|\}\subset\R^8
\]
is area-minimizing, although it is nonflat and has an isolated singularity at the origin. 

More recently, Savin established the corresponding rigidity theorem for Allen--Cahn minimizers~\cite{Savin2009}: if $n\leq7$, every nonconstant bounded global minimizer of the Allen--Cahn energy $\E$ is one-dimensional, i.e. it has the form
\[
u(x)=\tanh\left(\frac{x\cdot\nu-a}{\sqrt2}\right)
\quad\text{for some }\nu\in\mathbb S^{n-1}\text{ and }a\in\R.
\]

Because of the close relation between the Allen--Cahn energy and area, the existence of nonflat area-minimizing cones in dimension eight suggests that the dimension restriction in Savin's rigidity theorem should be sharp. More specifically, it is natural to try and construct non-one-dimensional Allen--Cahn solutions whose nodal sets are modelled on nonflat area-minimizing cones. In particular, for $m\geq 4$ consider the Simons cone
\begin{equation}\label{eq:Simons-cone-2m-1}
\mathcal C^{2m-1}=\{(x',x'')\in\R^m\times\R^m:|x'|=|x''|\}\subset \R^{2m}
\end{equation}
and the problem of constructing a bounded entire solution of \eqref{eq:AC} vanishing exactly on this cone, positive on $\{|x'|>|x''|\}$ and negative on $\{|x'|<|x''|\}$. The symmetries of the cone suggest seeking a solution depending only on $r=|x'|$ and $s=|x''|$, and changing sign under interchange of $r$ and $s$.

In~\cite{CabreTerra2009}, Cabr\'e and Terra proved that such solutions indeed exist, in fact in every even dimension $2m\geq2$. Their construction proceeds as follows. First, they minimize the Allen–Cahn energy on $\Omega_R=B_R\cap\{r>s\}$ among functions in $H^1_0(\Omega_R)$ that depend only on $r$ and $s$, imposing zero boundary values both on the cone and on the outer spherical boundary. Choosing a nonnegative minimizer $u_R$, they extend it to $B_R$ by setting $u_R(r,s)=-u_R(s,r)$ for $r<s$. They then show that, along a sequence $R\to\infty$, these extensions converge locally in $C^2$ to a nontrivial entire solution $U_{\mathrm{sad}}$, which satisfies $U_{\mathrm{sad}}(r,s)=-U_{\mathrm{sad}}(s,r)$ and has precisely the prescribed signs. Solutions with these properties are called \emph{saddle-shaped}.
Furthermore, in~\cite{CabreTerra2009,CabreTerra2010}, Cabr\'e and Terra studied the asymptotic behaviour of saddle-shaped solutions and their instability in dimensions four and six, and, finally, in~\cite{Cabre2012}, Cabr\'e proved uniqueness of $U_{\mathrm{sad}}$ within the class of saddle solutions, and its stability for $2m\geq14$. Further pointwise estimates were then obtained in~\cite{LWWPointwise2025}. 

The area-minimizing property of the cone $\Cone^{2m-1}$ for $m\geq4$ suggests that the saddle solution should be a global minimizer in the same dimensional range. This, however, does not follow from~\cite{CabreTerra2009,CabreTerra2010,Cabre2012,LWWPointwise2025}. Indeed, the variational construction a priori imposes rotational symmetry and fixes the solution to vanish on the cone, so it does not prove minimality under arbitrary compactly supported perturbations, including those that alter the nodal set. 

A further geometric construction suggests an additional way to understand the minimizing property. In order to show that the 7-dimensional Simons cone $\Cone^7$ is area-minimizing, Bombieri, De Giorgi, and Giusti~\cite{BDG1969} constructed smooth area-minimizing hypersurfaces $\Sigma_\pm$ on the two sides of $\Cone^7$ whose dilations $\{\lambda \Sigma_\pm\}_{\lambda>0}$ foliate $\R^8\setminus\Cone^7$. As $\lambda$ tends to zero, the leaves $\lambda\Sigma_\pm$ converge to the cone. The foliation $\{\lambda \Sigma_\pm\}_{\lambda>0}$ defines consistently oriented unit normal field on $\R^8\setminus\Cone^7$; this vector field is divergence-free and provides a calibration for $\Cone^7$. Later, Hardt and Simon extended this picture to all area-minimizing cones with an isolated singularity~\cite{HardtSimon1985}. 
An Allen--Cahn analogue would be a family $\{u_\lambda\}_{\lambda\in\R}$ of entire solutions of \eqref{eq:AC}, depending continuously on $\lambda$ in $C^2_{\mathrm{loc}}(\R^{2m})$, with $u_0=U_{\mathrm{sad}}$ and
\[
\begin{gathered}
	u_\lambda(x)<u_\mu(x)
	\quad\text{for every }x\in\R^{2m}\text{ whenever }\lambda<\mu,\\
	u_\lambda\to\pm1
	\quad\text{locally uniformly as }\lambda\to\pm\infty.
\end{gathered}
\]
Notice that these conditions ensure that every point of $\R^{2m}\times(-1,1)$ lies on the graph of exactly one member of the family. Once such a family is constructed, a sliding comparison with energy minimizers on bounded domains would then show that each $u_\lambda$ is a global minimizer~\cite{AlbertiAmbrosioCabre2001,Wang2018Slides}. Unlike the minimal hypersurface equation, \eqref{eq:AC} is not invariant under spatial dilations: rescaling a solution changes the equation's scale. Constructing such a family therefore requires more than dilating a single solution.

In~\cite{PW}, Pacard and Wei constructed stable Allen--Cahn solutions near sufficiently large dilations of the smooth area-minimizing leaves $\Sigma_\pm$ of the Hardt--Simon foliation associated with the Simons cone $\Cone^{2m-1}\subset\R^{2m}$ for $m\geq 4$. Liu, Wang, and Wei established the strict pointwise ordering of the two families of stable solutions constructed by Pacard and Wei, and used suitable members as upper and lower barriers in a minimization argument on expanding bounded domains to construct global minimizers whose nodal sets are asymptotic to suitable strictly area-minimizing Lawson cones, including the Simons cone in $\R^8$~\cite{LWW}. Their result establishes the sharpness of the dimension bound in Savin's theorem. For the Simons cones, Liu, Wang, and Wei construct, for each $a\in\R$, a global minimizer $U_a$ whose nodal set, outside a compact set, is described in the radial variables $r=|x'|$ and $s=|x''|$ by
\[
s=F_a(r),
\qquad
F_a(r)=r+a r^{-\alpha}+o(r^{-\alpha})
\quad\text{as }r\to\infty,
\]
where $\alpha$ is the decay exponent defined in \eqref{eq:exponents}; in particular, $\alpha=2$ in dimension eight. The coefficient $a$ is prescribed through the boundary data on expanding domains: these data approximate a transition across a geometric leaf whose asymptotic graph $s=f_a(r)$ has the same coefficient, with the cone itself used when $a=0$. Their asymptotic estimates show that the limiting nodal graph satisfies
\[
F_a(r)-f_a(r)=o(r^{-\alpha}),
\]
so the prescribed coefficient is retained in the entire limit~\cite[Propositions~7 and~13]{LWW}.

For $a=0$, this gives $F_0(r)-r=o(r^{-\alpha})$, which specifies a decay rate without asserting that $F_0(r)=r$. Proving that every global minimizer with this asymptotic behaviour and the prescribed phase orientation is the saddle-shaped solution therefore requires a further rigidity argument. Moreover, existence for each $a$ does not establish that the minimizers can be chosen so that $a\mapsto U_a$ is continuous in $C^2_{\mathrm{loc}}(\R^{2m})$.

In this paper, we resolve the rigidity question at $a=0$. For every $m\geq4$, we prove that the class $\mathcal M_0$ of global minimizers with zero asymptotic coefficient and the prescribed phase orientation consists of exactly one solution, and identify it with the saddle-shaped solution; see Definition~\ref{def:moduli} and Theorem~\ref{thm:main}. In particular, the saddle-shaped solution is a global minimizer in every even dimension $2m\geq8$. Our result thus identifies the unique member at the central parameter of the proposed family, without requiring a continuous choice of minimizers for all $a$. We now introduce the notation and precise asymptotic conditions needed to state the theorem.

\subsection*{Setup and main results} To state the rigidity result precisely, we describe the nodal set outside a compact set as a normal graph over the Simons cone and specify the asymptotic behaviour of its graphing function. Fix an integer $m\ge4$, set $n=2m$, and write $\Cone=\Cone^{2m-1}$ for the cone in \eqref{eq:Simons-cone-2m-1}. For $x=(x',x'')\in\R^m\times\R^m$, let
\[
r=|x'|,\qquad s=|x''|,\qquad
\rho=\frac{r+s}{\sqrt2},\qquad t=\frac{r-s}{\sqrt2},\qquad R=|x|.
\]
Here $t$ is the signed distance to $\Cone$, positive on $\{r>s\}$, and $\rho$ is the distance from the origin to a nearest point on the cone. In particular,
\begin{equation}\label{eq:geometry}
R^2=\rho^2+t^2,\qquad |t|\le\rho,\qquad
R/\sqrt2\le\rho\le R,\qquad \dist(x,\Cone)=|t|.
\end{equation}
The functions $\rho$ and $t$ are globally $1$-Lipschitz. See Appendix~\ref{app:identities} for a number of useful differential identities they satisfy.
For $|t|<\rho$ and $\theta_1,\theta_2\in\Sph^{m-1}$, we have
\[
x=\left(\frac{\rho+t}{\sqrt2}\theta_1,
\frac{\rho-t}{\sqrt2}\theta_2\right).
\]
For fixed $\rho, \theta_1,$ and $\theta_2$, varying $t$ yields a unit-speed line normal to the cone.

We now specify the exponent $\alpha$ appearing in the asymptotic expansion $F_a(r)=r+ar^{-\alpha}+o(r^{-\alpha})$ discussed above. To explain its geometric origin, consider the Jacobi equation $(\Delta_\Cone+|A_\Cone|^2)\xi=0$ on the cone $\Cone$. If $\xi$ is assumed to be invariant under $O(m)\times O(m)$, it depends only on the radial coordinate $\rho=|y|$, and the equation becomes
\[
\xi''+\frac{n-2}{\rho}\xi'+\frac{n-2}{\rho^2}\xi=0.
\]
Seeking solutions of the form $\xi(\rho)=\rho^{-\gamma}$ gives $\gamma^2-(n-3)\gamma+(n-2)=0$. For $n=2m\ge8$, this polynomial has two distinct positive roots $\alpha<\beta$. We set
\begin{equation}\label{eq:exponents}
\alpha=\frac{n-3-\sqrt{(n-3)^2-4(n-2)}}2,
\qquad \beta=n-3-\alpha,
\qquad c_m=2^{(\alpha-1)/2}.
\end{equation}
Thus $\rho^{-\alpha}$ and $\rho^{-\beta}$ are the two radial solutions of the Jacobi equation, with $\rho^{-\alpha}$ being the slowly decaying one. 
The exponents $\alpha, \beta$ satisfy $1<\alpha\le2$, $\alpha+\beta=n-3$, and $\alpha\beta=n-2$; for $m=4$, we have $(\alpha,\beta)=(2,3)$ and $c_m=\sqrt2$. 

Hereafter, all constants may depend on $m$. The dimensional restriction $m\ge4$ is essential to our argument, as it ensures both the strict minimality of the cone and the existence of two distinct real decay exponents. No assertions are made for $m\le3$.

\medskip

We shall now introduce the moduli space under consideration and present our main theorem.

\begin{definition}\label{def:moduli}
	For each $a\in\R$, let $\M_a$ denote the set of global minimizers
	$u:\R^n\to(-1,1)$ possessing the following properties: in the exterior region
	$\{\rho > R_u\}$, the entire nodal set is a single continuous graph
	$t=f_u(\rho,\theta_1,\theta_2)$ satisfying $|f_u|<\rho/2$, and
	\begin{equation}\label{eq:phase}
		\sgn u(x)
		=\sgn(t-f_u(\rho,\theta_1,\theta_2)).
	\end{equation}
	Furthermore, we require the uniform asymptotic behavior
	\begin{equation}\label{eq:coefficient}
		\sup_{\theta_1,\theta_2}\rho^\alpha
		\left|f_u(\rho,\theta_1,\theta_2)+c_m a\rho^{-\alpha}\right|
		\to 0\quad \text{as }\rho\to\infty.
	\end{equation}
\end{definition}

Note that for a symmetric graph $s=r+ar^{-\alpha}+o(r^{-\alpha})$, the corresponding normal height is exactly $f_u=-c_m a\rho^{-\alpha}+o(\rho^{-\alpha})$. Our main theorem is the following.

\begin{theorem}
\label{thm:main}
The set $\M_0$ consists of exactly one global minimizer $U_0$. This minimizer is
$O(m)\times O(m)$-invariant and satisfies
\[
	\begin{gathered}
	U_0(r,s)=-U_0(s,r),\\
	U_0(r,s)>0\quad\text{for }r>s,\qquad
	U_0(r,s)<0\quad\text{for }r<s.
	\end{gathered}
\]
In particular, $U_0$ is a saddle-shaped global minimizer.
\end{theorem}

Combining Theorem~\ref{thm:main} with the uniqueness of the saddle-shaped solution~\cite[Theorem~1.2]{Cabre2012}, we obtain the following consequence.

\begin{corollary}
	The saddle-shaped solution of \cite{CabreTerra2009} is a global minimizer.
\end{corollary}

In particular, our result establishes the uniqueness of the fast mode coefficient among minimizers, whenever the slow mode coefficient for the normal height is zero.

\subsection*{Method of the proof}
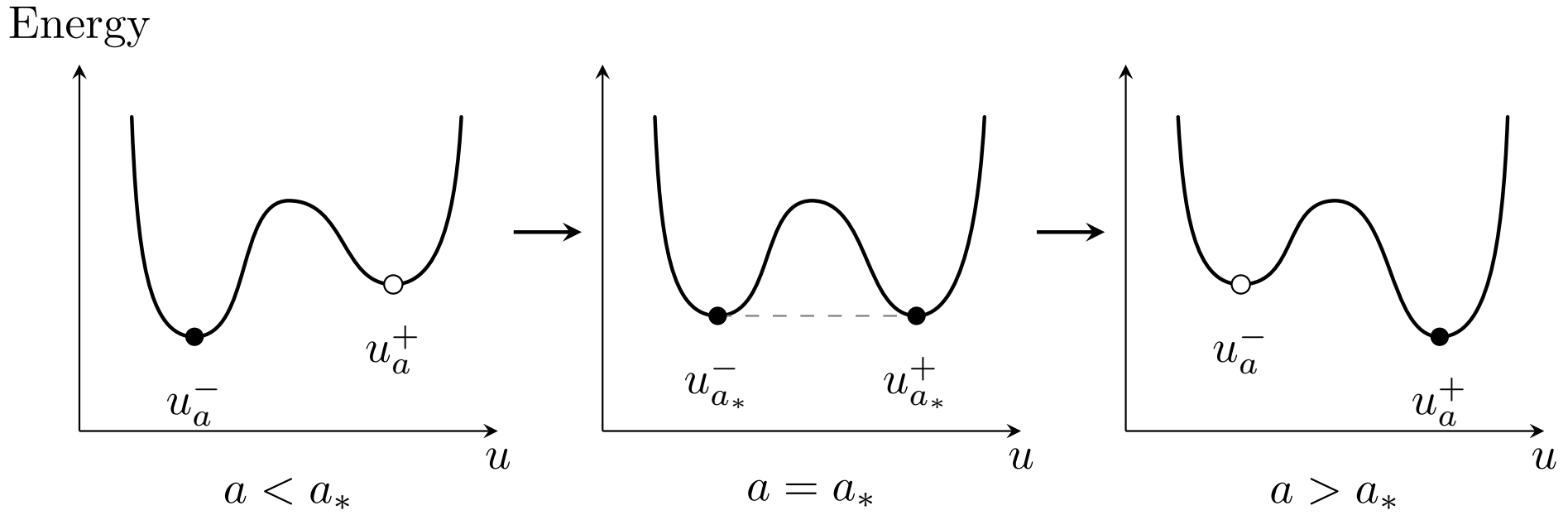
\begin{figure}[htbp]
	\centering
	\resizebox{\linewidth}{!}{%
		\begin{tikzpicture}[
			x=0.8cm,y=0.8cm,
			>=stealth,
			every node/.style={font=\small},
			curve/.style={thick,line cap=round}
			]
			
			\draw[->] (0,0) -- (4,0) node[below] {$u$};
			\draw[->] (0,0) -- (0,3.5) node[above] {Energy};
			\node at (2,-0.6) {$a<a_*$};
			
			\draw[curve]
			(0.5,3)
			.. controls (0.55,1.7) and (0.65,0.9) .. (1.1,0.9)
			.. controls (1.65,0.9) and (1.5,2.2) .. (2,2.2)
			.. controls (2.55,2.2) and (2.5,1.4) .. (3,1.4)
			.. controls (3.45,1.4) and (3.6,2.1) .. (3.65,3);
			
			\fill (1.1,0.9) circle (2pt);
			\draw[fill=white] (3,1.4) circle (2pt);
			\node[below=5pt] at (1.1,0.9) {$u_a^-$};
			\node[below=5pt] at (3,1.4) {$u_a^+$};
			
			\begin{scope}[shift={(5,0)}]
				\draw[->] (0,0) -- (4,0) node[below] {$u$};
				\draw[->] (0,0) -- (0,3.5) node[above] {};
				\node at (2,-0.6) {$a=a_*$};
				
				\draw[curve]
				(0.5,3)
				.. controls (0.55,1.8) and (0.65,1.1) .. (1.1,1.1)
				.. controls (1.65,1.1) and (1.5,2.2) .. (2,2.2)
				.. controls (2.55,2.2) and (2.5,1.1) .. (3,1.1)
				.. controls (3.45,1.1) and (3.6,2.1) .. (3.65,3);
				
				\draw[dashed,gray] (1.1,1.1) -- (3,1.1);
				\fill (1.1,1.1) circle (2pt);
				\fill (3,1.1) circle (2pt);
				\node[below=5pt] at (1.1,1.1) {$u_{a_*}^-$};
				\node[below=5pt] at (3,1.1) {$u_{a_*}^+$};
			\end{scope}
			
			\begin{scope}[shift={(10,0)}]
				\draw[->] (0,0) -- (4,0) node[below] {$u$};
				\draw[->] (0,0) -- (0,3.5) node[above] {};
				\node at (2,-0.6) {$a>a_*$};
				
				\draw[curve]
				(0.5,3)
				.. controls (0.55,2.1) and (0.65,1.4) .. (1.1,1.4)
				.. controls (1.65,1.4) and (1.5,2.2) .. (2,2.2)
				.. controls (2.55,2.2) and (2.5,0.9) .. (3,0.9)
				.. controls (3.45,0.9) and (3.6,1.7) .. (3.65,3);
				
				\draw[fill=white] (1.1,1.4) circle (2pt);
				\fill (3,0.9) circle (2pt);
				\node[below=5pt] at (1.1,1.4) {$u_a^-$};
				\node[below=5pt] at (3,0.9) {$u_a^+$};
			\end{scope}
			
			\draw[->,thick] (4.15,1.9) -- (4.8,1.9);
			\draw[->,thick] (9.15,1.9) -- (9.8,1.9);
			
		\end{tikzpicture}%
	}
	\caption{Schematic illustration of a possible jump in the global minimizers. Two branches of local minimizers may persist continuously while their energies cross at $a=a_*$. At the crossing, two distinct global minimizers coexist.}
	\label{fig:minimizer-switching}
\end{figure}

Our starting point is the observation that the global minimizers constructed in \cite{LWW} are ordered, as discussed in \cite[p.~36]{Wang2018Slides}. We rely on the existence of reference minimizers with a non-zero nodal coefficient, as established in \cite[Theorem~1]{LWW} (see Section~\ref{sec:ordering}). Building on this, we prove rotational symmetry for every member of $\M_a$, alongside strict ordering between different parameters:
\[
a<b,\quad u\in\M_a,\quad v\in\M_b
\quad\Longrightarrow\quad u<v
\quad\text{in }\R^n;
\]
see Lemma~\ref{lem:rotational-symmetry} and Proposition~\ref{prop:ordering}. The comparison argument combines asymptotic control with an exchange of the minimum and maximum of two minimizers. Ordering already implies that only countably many jumps can occur along a family of solutions $\{U_a\}_{a\in\R}$. To obtain a complete continuous family as above, the main difficulty is to rule out jumps at every parameter. Even when global minimizers are unique near $a_*$ ($|\M_a|=1$ for nearby $a\ne a_*$), this does not itself imply uniqueness in $\M_{a_*}$, as illustrated in Figure~\ref{fig:minimizer-switching}.

Thus, if the space $\M_a$ contains two distinct minimizers, the ordering between parameters precludes a continuous choice of minimizers through that parameter. One possible approach to rule out this phenomenon---and thereby exclude multiple minimizers within a given class---is to estimate the relative energy between two candidate minimizers and directly prove uniqueness. However, this approach is difficult due to the nonlinearity of the equation and the geometry of the minimizers and their nodal sets.

In the $O(m)\times O(m)$-symmetric setting with $a=0$, we can exploit reflection symmetry across the line $r=s$ by defining the operator
\[
(\mathcal R U)(r,s):=-U(s,r).
\]
The ordering property of minimizers yields an upper envelope $U_0^+$ and a lower envelope $U_0^-$ in $\M_0$, obtained as the canonical one-sided limits
\[
U_0^-=\lim_{a\uparrow0}U_a,
\qquad U_0^+=\lim_{a\downarrow0}U_a,
\qquad U_0^-\le V\le U_0^+\quad(V\in\M_0);
\]
see Lemma~\ref{lem:envelope}. Since $\mathcal R$ preserves $\M_0$ while reversing the pointwise order, it follows that
\[
\mathcal R U_0^+=U_0^-,
\qquad U_0^+-U_0^-\ge0;
\]
see Proposition~\ref{prop:reflection}. Thus equality of the two envelopes simultaneously gives uniqueness in $\M_0$ and the saddle antisymmetry.

The idea is to apply maximum principles---generalizing the qualitative properties of solutions in \cite{CabreTerra2010,Cabre2012}---to extract pointwise estimates from the additional condition $U_0^+-U_0^-\ge0$. A key step is to rewrite the Allen--Cahn equation as a coupled system in terms of the odd and even parts of $U=U_0^+$:
\[
U^*(r,s)=U(s,r),\qquad
S=\frac{U-U^*}{2},\qquad
T=\frac{U+U^*}{2}=\frac{U_0^+-U_0^-}{2}\ge0.
\]
The function $T$ measures the failure of antisymmetry. If $T\not\equiv0$, the strong maximum principle gives $T>0$. The pointwise estimates in Proposition~\ref{prop:pointwise-estimates} then yield
\[
S,S_\rho,S_t>0,\qquad T_\rho,T_t<0
\quad\text{on }0<t<\rho,
\]
and ensure that the relative energy density
\[
e_{\mathrm{rel}}=\frac{e(U_0^-)-e(U_0^+)}2
=-S_\rho T_\rho-S_tT_t+ST(1-S^2-T^2)
\]
is strictly positive on this side of the cone, where $e(u)=\frac12|\nabla u|^2+W(u)$. On reflection-invariant domains the total energy difference cancels, so we use a weighted quantity on one side.

The positivity of $e_{\mathrm{rel}}$ allows us to estimate the energy difference between $U_0^+$ and $U_0^-$. Our approach is partially inspired by the work of the first author and F.~Wang \cite{LeeWangPE} on the rigidity of conformally Einstein manifolds. There, second-order differential inequalities govern the propagation of integral curvature quantities and yield rigidity; see \cite[Section 5]{LeeWangPE}. For the Allen--Cahn equation, the variational structure provides additional control that strengthens this ODE argument, using the stress--energy tensor computed in \cite{Alikakos2011}.

We exploit the variational structure through the ``relative'' stress--energy tensor
\[
\mathsf T_{\mathrm{rel}}
=\frac12\bigl(\mathsf T[U_0^-]-\mathsf T[U_0^+]\bigr),
\qquad
\mathsf T[u]=\nabla u\otimes\nabla u-e(u)I_n.
\]
Since both functions solve the Allen--Cahn equation, $\mathsf T_{\mathrm{rel}}$ is divergence-free, and its conservation laws describe how the relative energy is distributed across $\R^n$; see Lemma~\ref{lem:stress}.

Let $\omega=(\rho^2-t^2)^{m-1}$ denote the volume density in $(\rho,t)$ coordinates, up to a dimensional constant. The ordinary relative energy
\[
\int_0^\rho\int_0^{\rho'}
\omega(\rho',t')e_{\mathrm{rel}}(\rho',t')\,dt'\,d\rho'
\]
is taken over a truncated portion of $\{r>s\}$, but does not directly give the differential inequality needed for rigidity. The conservation laws instead lead to the weighted quantity
\[
P(\rho)
=-2\int_0^\rho\int_0^{\rho'}
\partial_{t'}\omega(\rho',t')\,
e_{\mathrm{rel}}(\rho',t')\,dt'\,d\rho';
\]
see the discussion in the proof of Lemma~\ref{lem:rigidity}. Its normalized primitive
\[
F(\rho)=\rho^{-3}\int_0^\rho P(\rho')\,d\rho'
\]
satisfies the second-order differential inequality
\[
F''-\frac{n-8}{\rho}F'-\frac{n-8}{\rho^2}F>0.
\]
The condition $n\ge8$ enters here: positivity and regularity at the origin force $F(\rho)\ge c\rho$ for large $\rho$.

To obtain a contradictory upper bound, we establish refined asymptotic estimates for $U_0^\pm$ in Corollary~\ref{cor:sharp-profile} that sharpen the asymptotics obtained in \cite{LWW}. Let $H$ denote the one-dimensional transition profile $H(t)=\tanh(t/\sqrt2)$, and let $h=H'$. Projection onto the translation mode $h$ separates normal displacement from the orthogonal component, which is controlled by the spectral gap of the one-dimensional linearized operator. After subtracting the universal curvature correction, the projected equation has the radial Jacobi operator as its leading part. For two $O(m)\times O(m)$-invariant solutions trapped between reference minimizers and having the same slow coefficient, the difference equation permits an iteration to the fast rate $\rho^{-\beta}$ (Lemma~\ref{lem:fast-difference}). In particular, with the necessary derivative bounds,
\[
T=O_\mu(\rho^{-\beta}),\qquad
S-H-(n-2)\rho^{-2}\eta=O_\mu(\rho^{-4}),
\]
where $(-\partial_t^2+3H^2-1)\eta=-tH'$. Isolating this fast mode yields $P(\rho)=O(\rho^\alpha)$ and hence $F(\rho)=O(\rho^{\alpha-2})$. Since $\alpha\le2$, this is incompatible with the lower bound. We conclude that $T\equiv0$ and that the two extremal minimizers coincide.

This proves rigidity at the reflection-invariant parameter without the need for uniqueness at every parameter or a complete foliation. By the comparison principle and the equality $U_0^-=U_0^+$, we also obtain local smooth convergence to the saddle solution as $a\to0$, regardless of the choice of $U_a\in\M_a$. Finally, global minimality is inherited from the extremal limits, and the rigidity argument determines their common value.

\subsection*{Outline of the paper}
 In Section~\ref{sec:prelim}, we establish comparison principles and obtain
preliminary estimates for the solution profile near the cone.
Section~\ref{sec:refined} further refines these estimates to isolate the
contribution of the slow mode. In Section~\ref{sec:ordering}, we prove
the rotational symmetry and strict ordering of minimizers with distinct
asymptotic parameters, allowing us to construct the extremal minimizers
$U_0^\pm$ as one-sided limits as $a\to0$.
Section~\ref{sec:signs} establishes pointwise inequalities for the odd
and even parts of $U_0^+$ alongside their derivatives. Finally, in Section~\ref{sec:main},
we combine these inequalities with the conservation law for the relative
stress--energy tensor to demonstrate that $U_0^+=U_0^-$, thereby completing the proof
of the main theorem. Appendix~\ref{app:identities} collects the coordinate
identities, projection formulas, and linear estimates utilized throughout the paper.

\subsection*{AI disclosure}
AI tools were utilized to generate the figure in the introduction, ensure grammatical consistency, and verify the correctness of the algebraic computations.

\subsection*{Notation}
We collect the most frequently used notation below.

\begin{itemize}
\item[-] $\alpha,\beta,c_m$: The slow and fast decay exponents and the normalization constant of the asymptotic profile.

\item[-] $r,s$; $\rho,t$: $r=|x'|$, $s=|x''|$; $\rho=(r+s)/\sqrt2$, $t=(r-s)/\sqrt2$.

\item[-] $R$, $\Cone$: $R=|x|=(\rho^2+t^2)^{1/2}$ and $\Cone=\{t=0\}$. 

\item[-] $\omega=(\rho^2-t^2)^{m-1}$ is the Euclidean volume density in
$(\rho,t)$ coordinates up to a dimensional constant.

\item[-] $W(z)$, $\E(u;\Omega)$: The potential $W(z)=(1-z^2)^2/4$ and the Allen--Cahn energy over the domain $\Omega$.

\item[-] $e(u)=\tfrac12|\nabla u|^2+W(u)$, $\delta_u:=2W(u)-|\nabla u|^2$

\item[-] $f_u$: The normal height of the nodal set of $u$ over $\Cone$ for large $\rho$. 

\item[-] $\M_a$: The class consisting of global minimizers with the prescribed phase orientation, such that $f_u=-c_m a\rho^{-\alpha}+o(\rho^{-\alpha})$.

\item[-] $\mathcal T_{\kappa,R}$: The tubular region $\{\rho>R,\ |t|<\kappa\rho\}$.

\item[-] $D^j$, $O_\mu(\rho^{-k})$: Euclidean derivatives of order $j$; the bound $C\rho^{-k}e^{-\mu|t|}$ in the tubular region.

\item[-] $H$, $h$, $\sigma$: $H(t)=\tanh(t/\sqrt2)$, $h=H'$, and $\sigma=\int_\R h^2\,dt$.

\item[-] $L_H$, $\eta$: $L_H=-\partial_t^2+3H^2-1$ and the odd solution $L_H\eta=-th$ from Lemma~\ref{lem:normal}.

\item[-] $\mathcal L_U$, $\mathcal L_H$: The linearized operators $-\Delta+3U^2-1$ and $-\Delta+3H(t)^2-1$.

\item[-] $\mathcal L_0$, $\Pi$: $\mathcal L_0=-\partial_{\rho\rho}-(n-2)\rho^{-1}\partial_\rho+L_H$; $\Pi$ is the orthogonal projection onto $h^\perp$.

\item[-] $J$: The radial Jacobi operator over the cone $\Cone$.

\item[-] $\epsilon_*$: A fixed positive exponent strictly smaller than $\min\{(\alpha-1)/4,(\beta-\alpha)/4,1/2\}$.

\item[-] $\theta(\rho,t), \theta_N$, $\chi_\epsilon$: Cutoff functions.

\item[-] $\mathcal K$: The vector field $\mathcal K=s\partial_r-r\partial_s=-t\partial_\rho+\rho\partial_t$. 

\item[-] $U^*$, $\mathcal R U$: Reflection operators defined as $U^*(x',x'')=U(x'',x')$ and $\mathcal R U=-U^*$.

\item[-] $S,T$: The odd and even parts $(U-U^*)/2$ and $(U+U^*)/2$.
\end{itemize}

\section{Preliminary estimates}\label{sec:prelim}

In this section, we establish comparison principles and estimates
for the profiles outside a sufficiently large compact set. The coordinate identities used throughout
the paper are collected in Appendix~\ref{app:identities}.

\begin{lemma}\label{lem:flat}
Let $w$ be a stable entire solution of \eqref{eq:AC} with $|w|\le1$.
\begin{enumerate}[label=\textup{(\roman*)}]
\item If $w\ge0$, then $w\equiv1$; if $w\le0$, then $w\equiv-1$.
\item If $w=0$ on $\{x_n=0\}$, with $w\ge0$ in $\{x_n>0\}$
and $w\le0$ in $\{x_n<0\}$, then $w(x)=H(x_n)$.
\end{enumerate}
\end{lemma}

\begin{proof}
The classification in the whole space is due to Aronson and Weinberger
\cite{AW1978}, and the half-space classification is due to Angenent
\cite{Angenent1985}; stability rules out the identically zero solution.
\end{proof}

\begin{lemma}\label{lem:defect}
Let $u:\R^n\to(-1,1)$ be a stable entire solution whose
exterior nodal set is a continuous normal graph over the cone $\Cone$ with the phase
orientation \eqref{eq:phase}, and suppose that its height tends uniformly
to zero at infinity. Then
\begin{equation}\label{eq:defect}
\delta_u(x):=\frac12(1-u(x)^2)^2-|\nabla u(x)|^2
\to 0\qquad\text{uniformly as }|x|\to\infty.
\end{equation}
\end{lemma}

\begin{proof}
We argue by taking limits of translated solutions, following
\cite[proof of Theorem~1.6]{CabreTerra2010}. Suppose toward a contradiction that the asserted
uniform convergence fails. Then there exist $\varepsilon_0>0$
and a sequence of points $x_j\in\R^n$ with $|x_j|\to\infty$ such that
$
|\delta_u(x_j)|\ge\varepsilon_0.
$
Up to a subsequence, either
$\dist(x_j,\Cone)\to\infty$ or $\dist(x_j,\Cone)$ remains bounded.

In the first case, the translated solutions $u(x_j+\cdot)$
converge along a subsequence in $C^1_{\mathrm{loc}}$ to a stable
entire solution of fixed sign. Lemma~\ref{lem:flat}\textup{(i)}
then implies that this limit is identically $\pm1$, which yields a contradiction.

In the second case, for each $j$ we translate a point on the cone nearest to $x_j$ to
the origin and rotate the corresponding unit normal to $e_n$. The translated and rotated cones
then converge smoothly on compact sets to $\{x_n=0\}$, while the
transformed points corresponding to $x_j$ remain bounded. The assumptions on the nodal height and the phase
orientation imply that any subsequential limit vanishes on
$\{x_n=0\}$ and has the prescribed signs on either side.
By Lemma~\ref{lem:flat}\textup{(ii)}, the limit must be $H(x_n)$,
whose discrepancy function vanishes identically. Consequently,
$C^1_{\mathrm{loc}}$ convergence again implies $\delta_u(x_j)\to0$,
a contradiction.

\end{proof}

We next record an elementary maximum principle on unbounded domains.
\begin{lemma}\label{lem:basic-mp}
Let $D\subset\R^n$ be an open set and consider the operator
\[
\mathcal A=-\Delta+\mathbf b\cdot\nabla+c,
\qquad |\mathbf b|\le B_0,\qquad c\ge c_0>0.
\]
If $v\in C^2(D)\cap C^0(\overline D)$ is bounded from above and satisfies
\[
\mathcal Av\le0\quad\text{in }D,\qquad
v\le0\quad\text{on }\partial D,
\]
then $v\le0$ in $D$.
\end{lemma}

\begin{proof}
Set $\Phi(x)=C+\sqrt{1+|x|^2}$. For $C>0$ sufficiently large, we have
\[
\mathcal A\Phi\ge c_0(C+1)-n-B_0>0.
\]
For any $\varepsilon>0$, a positive maximum of
$v-\varepsilon\Phi$ over $\overline D$ must be attained at an interior point of $D$, since $v\le0$ on $\partial D$, $v$ is
bounded from above, and $\Phi(x)\to\infty$ as $|x|\to\infty$.
At such a point we have
$\mathcal A(v-\varepsilon\Phi)>0$, which contradicts
$\mathcal A(v-\varepsilon\Phi)<0$ in $D$.
Thus $v\le\varepsilon\Phi$ in $D$, and letting $\varepsilon\downarrow0$
yields the claim.
\end{proof}

We next construct a supersolution from functions of a known or approximate solution. For an exact solution $w$, we set
$\mathcal L_w=-\Delta+3w^2-1$ and denote its discrepancy function by
$\delta_w$, as in \eqref{eq:defect}. Direct calculation yields
\[
\mathcal L_w(1-w^2)=-2\delta_w,\qquad
\mathcal L_w(3w-1)=1-3w^2+6w^3\ge\frac12
\quad\text{for }0\le w\le 1.
\]
Thus $1-w^2$ lies approximately in the kernel of $\mathcal L_w$ when the defect is small,
whereas $\mathcal L_w(3w-1)$ is bounded away from zero. Adding twice the first
function to the second yields a function bounded below by $1$ on $[0,1]$:
\[
\Psi(w)=2(1-w^2)+3w-1=1+3w-2w^2\ge1,
\]
for which both $\Psi(w)$ and $\mathcal L_w \Psi(w)$ are positive provided $\delta_w$ is sufficiently small. In the following lemma, we
also accommodate a small error, which enables
us to construct barriers from $\pm H(t)$ in $(\rho, t)$-coordinates.

\begin{lemma}\label{lem:Psi-MP}
Let $D\subset\R^n$ be an open set, let $Q\in C^0(D)$, and let
$\mathcal L=-\Delta+Q$. Suppose that
$w\in C^2(D)\cap C^0(\overline D)$ satisfies $0\le w\le 1$,
$Q\ge3w^2-1$, and
\begin{equation}\label{eq:Psi-smallness}
4\big|\frac12(1-w^2)^2-|\nabla w|^2\big|
+3|-\Delta w-w+w^3|\le\frac14.
\end{equation}
If we define
$\Psi(s)=1+3s-2s^2$ for $0\le s\le 1$, then
\begin{equation}\label{eq:positive-psi}
1\le\Psi(w)\le3,\qquad
\mathcal L\Psi(w)\ge\frac14,
\end{equation}
and the conjugate operator defined by
\begin{equation}\label{eq:Psi-conjugation}
	\begin{gathered}
	\mathcal A f:=\frac{\mathcal L(\Psi(w)f)}{\Psi(w)}
	=-\Delta f+\mathbf b\cdot\nabla f+cf,\\
	\mathbf b=-2\nabla\log\Psi(w),\qquad
	c=\frac{\mathcal L\Psi(w)}{\Psi(w)},
	\end{gathered}
\end{equation}
satisfies $|\mathbf b|\le6$ and $c\ge\frac{1}{12}$.
Consequently, if $v\in C^2(D)\cap C^0(\overline D)$ is bounded from above,
$\mathcal Lv\le0$ in $D$, and $v\le0$ on $\partial D$, then $v\le0$ in $D$.
\end{lemma}

\begin{proof}
For $0\le s\le1$, a direct calculation yields
\begin{equation}\label{eq:Psi-bounds}
1\le\Psi(s)\le3,\qquad |\Psi'(s)|\le3,\qquad
1-3s^2+6s^3\ge\frac12.
\end{equation}
Applying the chain rule, we obtain
\begin{equation}\label{eq:P-identity}
	\begin{aligned}
	\mathcal L\Psi(w)
	={}&1-3w^2+6w^3
	-4\left(\frac12(1-w^2)^2-|\nabla w|^2\right)\\
	&+\Psi'(w)(-\Delta w-w+w^3)
	+(Q-3w^2+1)\Psi(w).
	\end{aligned}
\end{equation}
Consequently, \eqref{eq:Psi-smallness} implies
$\mathcal L\Psi(w)\ge1/4$, while the same hypothesis ensures that
$|\nabla w|^2\le1/2+1/16<1$.
Direct expansion of $\mathcal L(\Psi(w)f)$ then verifies
\eqref{eq:Psi-conjugation} and its coefficient bounds.
Finally, the quotient $v/\Psi(w)$ is bounded from above, is nonpositive on $\partial D$,
and satisfies $\mathcal A(v/\Psi(w))\le0$.
The desired conclusion now follows directly from Lemma~\ref{lem:basic-mp}.
\end{proof}

We next establish rough estimates both near and away from the cone.

\begin{proposition}\label{prop:away}
Let $u\in\M_a$. There exist constants $t_0, C, \nu > 0$ such that
\begin{equation}\label{eq:away}
|u(x)-\operatorname{sgn}(t)|\le Ce^{-\nu|t|}\qquad\text{for }|t|\ge t_0.
\end{equation}
\end{proposition}

\begin{proof}
Since $u$ is stable and its nodal heights are bounded,
the compactness of translated solutions together with Lemma~\ref{lem:flat}\textup{(i)}
implies that $u\to\pm1$ uniformly as $t\to\pm\infty$.
For sufficiently large positive $t$,
\[
\Delta(1-u)=u(1+u)(1-u)\ge1-u.
\]
Applying \cite[Lemma~B.1]{Wang2017} in balls of radius comparable
to $|t|$, and arguing analogously for $-u$, proves
\eqref{eq:away}.
\end{proof}

\begin{proposition}\label{prop:coarse}
Let $0<\kappa<1/4$. If $u\in\M_a$, there exist constants
$R,C,\mu>0$ depending on $u$ and $\kappa$ such that
\begin{equation}\label{eq:coarse}
\sum_{j=0}^2 |D^j(u-H(t))|
\le C\rho^{-\alpha}e^{-\mu|t|}
\qquad\text{in }\mathcal T_{\kappa,R}.
\end{equation}
\end{proposition}

\begin{proof}
By Proposition~\ref{prop:away} and the exponential decay of
$H(t)-\operatorname{sgn}(t)$, there exist $t_0>0$ and $0<\nu<\sqrt2$ such that
\[
|u-H(t)|\le Ce^{-\nu|t|}\qquad\text{for }|t|\ge t_0.
\]
On any fixed strip $|t|\le T$, the argument
of Lemma~\ref{lem:defect}, combined with
Lemma~\ref{lem:flat}\textup{(ii)}, implies that $u-H(t)\to0$
uniformly as $\rho\to\infty$. Combining this fact with the preceding
bound yields
\[
\sup_{\mathcal T_{2\kappa,R}}|u-H(t)|\to 0
\qquad\text{as }R\to\infty.
\]

We first work in the region $\mathcal T_{2\kappa,R}\cap\{t>0\}$. By
\eqref{apeq:gradient-identities} in Appendix~\ref{app:identities},
the gradients of $\rho$ and $t$ are orthonormal, with Laplacians given by
\[
\Delta\rho=A=\frac{(n-2)\rho}{\rho^2-t^2}
\quad\text{and}\quad
\Delta t=B=-\frac{(n-2)t}{\rho^2-t^2}.
\]
In particular, $|A|+|B|\le C/\rho$ throughout this region. Setting
\[
\mathcal L=-\Delta+u^2+uH+H^2-1,
\]
we find that $\mathcal L(u-H)=BH'$ since $\Delta H=H^3-H+BH'$, where
\begin{equation}\label{eq:forcing}
|BH'|\le C\rho^{-2}t e^{-\sqrt2 t}.
\end{equation}

The profile $H(t)$ has zero discrepancy since $|\nabla t|=1$,
while its residual error in the Allen--Cahn equation is $-BH'=O(\rho^{-2})$.
Consequently, Lemma~\ref{lem:Psi-MP} gives
\[
2|\nabla\log\Psi(H)|\le6,\qquad
\frac{(-\Delta+3H^2-1)\Psi(H)}{\Psi(H)}\ge\frac1{12}.
\]
Since $(u-H)(u+2H)=o(1)$ uniformly as $R\to\infty$, we find that
\[
\frac{\mathcal L\Psi(H)}{\Psi(H)}
=\frac{(-\Delta+3H^2-1)\Psi(H)}{\Psi(H)}
+(u-H)(u+2H)\ge\frac1{24}
\]
for all sufficiently large $R$. If $\mathcal Lv\le0$, then $f=v/\Psi(H)$ satisfies
\[
-\Delta f-2\nabla\log\Psi(H)\cdot\nabla f
+\frac{\mathcal L\Psi(H)}{\Psi(H)}f\le0.
\]
Recalling that $1\le\Psi(H)\le3$, Lemma~\ref{lem:basic-mp} applied to $f$
yields the desired maximum principle for $\mathcal L$ on functions bounded above.

Choose $0<\mu<\nu$ so that $\mu^2+6\mu\le1/96$, and put
$\varphi=\rho^{-\alpha}e^{-\mu t}$. Then
\[
\frac{\Delta\varphi}{\varphi}
=\mu^2-\mu B-\frac{\alpha A}{\rho}
+\frac{\alpha(\alpha+1)}{\rho^2}
=\mu^2+O(\rho^{-1}).
\]
Since $\Psi(H)$ depends only on $t$, the product rule gives
\begin{equation}\label{eq:barrier-calculation}
	\begin{aligned}
	\frac{\mathcal L(\varphi\Psi(H))}{\varphi\Psi(H)}
	&=\frac{\mathcal L\Psi(H)}{\Psi(H)}
	-\frac{\Delta\varphi}{\varphi}+2\mu\,\partial_t\log\Psi(H)\\
	&\ge\frac1{24}-\mu^2-6\mu-\frac{C}{\rho}
	\ge\frac1{48}
	\end{aligned}
\end{equation}
after increasing $R$. Moreover, by \eqref{eq:forcing},
\[
\frac{|Bh|}{\varphi\Psi(H)}
\le C\rho^{\alpha-2}t e^{-(\sqrt2-\mu)t}\le C,
\]
since $\alpha\le2$ and $\mu<\sqrt2$. Hence we may choose $C$
so that $\mathcal L(C\varphi\Psi(H))\ge|Bh|$.

We now compare $|u-H|$ with $C\varphi\Psi(H)$ on the boundary of $\mathcal T_{2\kappa,R}\cap\{t>0\}$.
On $\{t=0\}$, we have $\varphi\Psi(H)=\rho^{-\alpha}$, while
\eqref{eq:coefficient} together with the uniform gradient bound yields
\[
|u(\rho,0,\omega,\theta)|
\le C|f_u(\rho,\omega,\theta)|\le C\rho^{-\alpha}.
\]
On $\{t=2\kappa\rho\}$, the decay estimate away from the cone gives
\[
\frac{|u-H|}{\varphi\Psi(H)}
\le C\rho^\alpha e^{-2\kappa(\nu-\mu)\rho}\le C.
\]
Finally, on $\{\rho=R\}$, the corresponding bound follows from $|u-H|\le2$ and
$\varphi\Psi(H)\ge R^{-\alpha}e^{-2\kappa\mu R}>0$.

Upon choosing $C>0$ sufficiently large, the functions
$v_\pm=\pm(u-H)-C\varphi\Psi(H)$ are nonpositive on the boundary,
bounded from above, and satisfy $\mathcal Lv_\pm\le0$.
Lemma~\ref{lem:Psi-MP} therefore yields
$|u-H|\le C\varphi\Psi(H)\le3C\varphi$ in $\mathcal T_{2\kappa,R}\cap\{t>0\}$.
Applying the same argument to
$\widetilde u(x',x'')=-u(x'',x')$ establishes the corresponding bound for $t<0$,
and hence
\begin{equation}\label{eq:weighted-value}
|u-H(t)|\le C\rho^{-\alpha}e^{-\mu|t|}
\qquad\text{in }\mathcal T_{2\kappa,R}.
\end{equation}

Finally, the explicit expression for $B$ yields
\[
|Bh|+|\nabla(Bh)|
\le C\rho^{-2}(1+|t|)e^{-\sqrt2|t|}.
\]
Since the weight $\rho^{-\alpha}e^{-\mu|t|}$ in \eqref{eq:weighted-value} varies by at most a bounded factor on
unit balls and the coefficients of $\mathcal L$ enjoy uniform
local H\"older bounds, interior Schauder estimates applied to
$\mathcal L(u-H)=Bh$ yield \eqref{eq:coarse}
in $\mathcal T_{\kappa,R}$ after increasing $R$ if necessary.
\end{proof}

\section{Refined profiles of \texorpdfstring{$O(m)\times O(m)$}{O(m)xO(m)}-invariant solutions}
\label{sec:refined}

The goal of this section is to establish refined asymptotic profiles for the global minimizers, particularly for $U_0^\pm$. To prove the uniqueness of minimizers, we must extract the fast decay exponent $\beta$ from the asymptotic expansion of these solutions.

For $\kappa\in(0,1/4)$, $R\ge1$, and $\mu>0$, we write
$E=O_\mu(\rho^{-k})$ in $\mathcal T_{\kappa,R}$ if
\[
\sup_{x\in\mathcal T_{\kappa,R}}
\rho(x)^k e^{\mu|t(x)|}|E(x)|<\infty,
\]
where $\rho$ and $t$ are viewed as functions on $\R^n$. For tensor-valued quantities, we use the
Euclidean norm and specify the orders of differentiation explicitly; for
$C^{0,\gamma}$ estimates on unit balls, we follow the convention introduced at the beginning of Appendix~\ref{app:identities}.
Fix
\[
0<\kappa<1/4,\qquad
0<\epsilon_*<\min\{(\alpha-1)/4,(\beta-\alpha)/4,1/2\}.
\]

In all statements below concerning entire Allen--Cahn solutions,
the $O_\mu$ estimates hold in $\mathcal T_{\kappa,R}$ for some
$\mu\in(0,\sqrt2)$ and all sufficiently large $R\ge1$; we leave this region implicit in what follows. These parameters, together with the threshold
$t_0\ge1$ in the estimates away from the cone and the implicit constants, depend only on
$n,\kappa,\epsilon_*$, and the bounding minimizers specified in
each statement, and are independent of the intermediate solutions under consideration.

The projection identities, orthogonal estimates, and localization lemmas
used below are collected in Appendix~\ref{app:identities}.

\begin{proposition}\label{prop:trapping-bounds}
Let $U_{b_\pm}\in\M_{b_\pm}$ be minimizers.
Every entire solution $u$ of \eqref{eq:AC} such that
$U_{b_-}\le u\le U_{b_+}$ in $\R^n$ satisfies
\[
	\begin{aligned}
	D^j(u-H)&=O_\mu(\rho^{-\alpha})
	&&\text{in }\mathcal T_{\kappa,R}\quad\text{for }0\le j\le2,\\
	u-\operatorname{sgn}(t)&=O(e^{-\mu|t|})
	&&\text{for }|t|\ge t_0.
	\end{aligned}
\]
\end{proposition}

\begin{proof}
Propositions~\ref{prop:away} and~\ref{prop:coarse}
give the $C^0$ bounds. The $C^2$ estimates follow from interior Schauder estimates. Note that the implied constants depend
only on the bounding minimizers.
\end{proof}

Let $J$ denote the radial part of the Jacobi operator on the Simons
cone $\Cone$. We will iteratively apply the following ODE estimate
to extract the decay rate $O_\mu(\rho^{-\beta})$.

\begin{lemma}\label{lem:euler}
Let $R\ge2$, $k>\alpha$, and let $\xi\in C^2([R,\infty))$ satisfy
\begin{equation}\label{eq:euler-operator}
J\xi:=\xi''+\frac{n-2}{\rho}\xi'
+\frac{n-2}{\rho^2}\xi=F(\rho),
\qquad |F(\rho)|\le M\rho^{-k-2}.
\end{equation}
Then the limit $A=\lim_{\rho\to\infty}\rho^\alpha\xi(\rho)$ exists.
Setting $E=\xi-A\rho^{-\alpha}$, we have, for $0\le j\le2$,
\begin{equation}\label{eq:euler-asymptotics}
\rho^j E^{(j)}=
	\begin{cases}
	O(\rho^{-k}) & \text{if } k<\beta,\\
	O(\rho^{-\beta}\log\rho) & \text{if } k=\beta,\\
	O(\rho^{-\beta}) & \text{if } k>\beta
	\end{cases}
\qquad\text{for }\rho\ge R.
\end{equation}
Furthermore, $|A|$ and the implicit constants in the $O(\cdot)$ bounds are
bounded by
\[
C_0(n,k,R)\bigl(M+|\xi(R)|+R|\xi'(R)|\bigr).
\]
\end{lemma}

\begin{proof}
The solutions to the associated homogeneous equation are $\rho^{-\alpha}$ and $\rho^{-\beta}$.
By the variation of constants formula \cite[Theorem~3.12]{Teschl2012},
\begin{equation}\label{eq:euler-tail}
	\begin{aligned}
	\xi-A\rho^{-\alpha}={}&C_2\rho^{-\beta}
	-\frac{\rho^{-\alpha}}{\beta-\alpha}
	\int_\rho^\infty s^{\alpha+1}F(s)\,ds\\
	&-\frac{\rho^{-\beta}}{\beta-\alpha}
	\int_R^\rho s^{\beta+1}F(s)\,ds.
	\end{aligned}
\end{equation}
Directly evaluating and bounding these two integrals establishes the existence of the limit defining $A$ together with the decay bounds on $E$. Differentiating \eqref{eq:euler-tail} then yields the corresponding bounds on the derivatives.
\end{proof}

\begin{lemma}\label{lem:coarse-expansion}
Let $U_{b_\pm}\in\M_{b_\pm}$ be minimizers, and let $u$ be an $O(m)\times O(m)$-invariant entire solution of \eqref{eq:AC}
satisfying $U_{b_-}\le u\le U_{b_+}$ in $\R^n$. Then
\begin{equation}\label{eq:refined-profile}
	\begin{gathered}
	u=H+A\rho^{-\alpha}h+(n-2)\rho^{-2}\eta+E,\\
	D^jE=O_\mu(\rho^{-\alpha-\epsilon_*})\qquad\text{for }0\le j\le2,
	\end{gathered}
\end{equation}
where the coefficient $A$ is given by
\[
A=\frac1{h(0)}\lim_{\rho\to\infty}\rho^\alpha u(\rho,0).
\]
In particular, if $u\in\M_a$, then $A=c_m a$.
\end{lemma}

\begin{proof}
Set $q=u-H$, and let $\theta(\rho,t)=\chi(t/\rho)$ denote the even
cutoff function from Lemma~\ref{lem:localization}.
Then $\widetilde q=\theta q$ coincides with $q$ on $\{|t|\le\kappa\rho\}$
and extends by zero to all $t\in\R$. For each $\rho$, we define
$
\xi(\rho)=\frac1\sigma\int_\R\tilde q(\rho,t)h(t)\,dt
$
and decompose $\widetilde q$ as
\[
\widetilde q=\xi h+(n-2)\rho^{-2}\eta+p.
\]
Since $\eta$ is odd and $h$ is even, it follows that $p(\rho,\cdot)\perp h$.
Proposition~\ref{prop:trapping-bounds} provides the coarse
estimates for $q$, which in turn imply that $|\xi|=O(\rho^{-\alpha})$.

Recall from \eqref{apeq:operators} that
$\mathcal L_0=-\partial_{\rho\rho}-\frac{n-2}{\rho}\partial_\rho+L_H$. Projecting $\mathcal L_0\widetilde q$ onto $h$ and using the fact that $L_Hh=0$,
we obtain
\[
J\xi=\frac{n-2}{\rho^2}\xi
-\frac1\sigma\int_\R(\mathcal L_0\widetilde q)h\,dt.
\]
Lemma~\ref{lem:coarse} then provides the following estimates for $p$, $\xi$, and $J\xi$:

\begin{equation*}
D^jp=O_\mu(\rho^{-2\alpha})\quad\text{for }0\le j\le2,
\qquad |\xi|=O(\rho^{-\alpha}),
\tag{\ref*{apeq:coarse-orthogonal}}
\end{equation*}
\begin{equation*}
J\xi=O(\rho^{-3\alpha}+\rho^{-\alpha-3})
=O(\rho^{-\alpha-2-2\epsilon_*}).
\tag{\ref*{apeq:coarse-Jacobi}}
\end{equation*}

Applying Lemma~\ref{lem:euler} with
$k=\alpha+2\epsilon_*<\beta$, we obtain a constant $A$ such that
\[
\sum_{j=0}^2\rho^j
\left|\frac{d^j}{d\rho^j}(\xi-A\rho^{-\alpha})\right|
\le C\rho^{-\alpha-2\epsilon_*}.
\]
On the region $\{|t|\le\kappa\rho\}$, where $\theta=1$, the remainder in
\eqref{eq:refined-profile} takes the form
\[
E=(\xi-A\rho^{-\alpha})h+p.
\]
Since $2\alpha>\alpha+2\epsilon_*$, the estimates on $\xi-A\rho^{-\alpha}$ combined with
the exponential decay of $h$ and its derivatives yield the asserted bounds for $E$.

At $t=0$, both $H$ and $\eta$ vanish, so
$\rho^\alpha u(\rho,0)\to Ah(0)$ as $\rho\to\infty$. Finally, if $u\in\M_a$, then
$f_u=O(\rho^{-\alpha})$ by \eqref{eq:coefficient}.
Evaluating the expansion at $t=f_u$ and using the asymptotic relations
$H(t)=h(0)t+o(t)$ and $\eta(t)=O(t)$ near $t=0$, we obtain
\[
0=u(\rho,f_u)
=h(0)\bigl(f_u+A\rho^{-\alpha}\bigr)+o(\rho^{-\alpha}).
\]
Comparing this with \eqref{eq:coefficient} yields $A=c_m a$.
\end{proof}

\begin{lemma}\label{lem:fast-difference}
Let $u$ and $v$ be $O(m)\times O(m)$-invariant entire solutions of \eqref{eq:AC} satisfying
\[
U_1\le u\le U_2,\qquad U_3\le v\le U_4
\quad\text{in }\R^n,
\]
where $U_i\in\M_{b_i}$ for $1\le i\le4$.
If
\[
\lim_{\rho\to\infty}\rho^\alpha
\bigl(u(\rho,0)-v(\rho,0)\bigr)=0,
\]
then the difference $w=u-v$ satisfies
\begin{equation}\label{eq:fast-difference}
	\begin{aligned}
	D^jw&=O_\mu(\rho^{-\beta})&&\text{for }0\le j\le2,\\
	w&=O(e^{-\mu|t|})&&\text{for }|t|\ge t_0.
	\end{aligned}
\end{equation}
\end{lemma}

\begin{proof}
By Lemma~\ref{lem:coarse-expansion}, we have $A_u=A_v$, so that
$w=u-v=O_{\mu_0}(\rho^{-\alpha-\epsilon_*})$.
Using the same cutoff function $\theta$, we set
\[
	\Lambda=(1-t^2/\rho^2)^{(n-2)/4},\qquad
	\widetilde w=\theta\Lambda w=\xi h+p,\qquad p\perp h.
\]
Multiplying by $\Lambda$ removes the first derivative in $t$ and
leaves $(n-2)\rho^{-1}\partial_\rho$ as the first-order radial term
in the Laplacian. Letting $Q_\Lambda$ denote the potential defined in
Lemma~\ref{lem:difference}, that lemma yields
\begin{equation*}
\mathcal L_0\widetilde w=-Q_\Lambda\widetilde w+\mathcal C_\chi[\Lambda w],
\qquad \|Q_\Lambda(\rho,\cdot)\|_{L^\infty}\le C\rho^{-\alpha},
\tag{\ref*{apeq:difference-equation}}
\end{equation*}
\begin{equation*}
|J\xi|\le C\rho^{-2-2\epsilon_*}|\xi|
+C\rho^{-\alpha}\|p\|_{L^2_t}+Ce^{-c\rho}.
\tag{\ref*{apeq:difference-Jacobi}}
\end{equation*}

Suppose that $w=O_\mu(\rho^{-k})$ for some $\alpha<k\le\beta$.
By the decay estimates away from the cone, this bound also holds on the support of the
cutoff function, possibly after decreasing $\mu$. Taking the component of \eqref{apeq:difference-equation}
orthogonal to $h$, we obtain
\[
\mathcal L_0p
=-\Pi(Q_\Lambda\widetilde w)+\Pi\mathcal C_\chi[\Lambda w]
=O_\mu(\rho^{-k-\alpha}).
\]
Thus, Lemma~\ref{lem:orthogonal} yields
\[
\|p(\rho,\cdot)\|_{L^2_t}
+\sup_t e^{\mu|t|}|p(\rho,t)|
\le C\rho^{-k-\alpha}.
\]
Consequently, since $2\alpha>2+2\epsilon_*$, it follows that
\[
|J\xi|\le C\rho^{-k-2-2\epsilon_*},\qquad
\lim_{\rho\to\infty}\rho^\alpha\xi=0.
\]

Lemma~\ref{lem:euler} improves the decay of $\xi$, and
hence of $w$, to $\rho^{-\min\{k+2\epsilon_*,\beta\}}$.
If $k+2\epsilon_*=\beta$, we can apply the iterate once more from $k+\epsilon_*$. Starting with $k=\alpha+\epsilon_*$ and
iterating finitely many times, we obtain
$w=O_\mu(\rho^{-\beta})$.

Interior Schauder estimate gives the $C^2$-estimate. 
The decay away from the cone follow from Proposition~\ref{prop:trapping-bounds}.
\end{proof}

In the following corollary, we show that after subtracting the first curvature correction $(n-2)\rho^{-2}\eta$ from the solution, the odd remainder is $O_\mu(\rho^{-4})$ and the even part is $O_\mu(\rho^{-\beta})$.

\begin{corollary}
\label{cor:sharp-profile}
Let $V\in\M_0$ be $O(m)\times O(m)$-invariant, satisfying
$U_{b_-}\le V\le U_{b_+}$ in $\R^n$, where
$U_{b_\pm}\in\M_{b_\pm}$ are minimizers. Define
\[
	\begin{gathered}
	V^*(\rho,t)=V(\rho,-t),\qquad
	S=(V-V^*)/2,\qquad T=(V+V^*)/2,\\
	E=S-H-(n-2)\rho^{-2}\eta.
	\end{gathered}
\]
Then
\begin{equation}\label{eq:zero-odd-refinement}
D^jE=O_\mu(\rho^{-4}),\qquad D^jT=O_\mu(\rho^{-\beta})
\qquad\text{for }0\le j\le2.
\end{equation}
\end{corollary}

\begin{proof}
Both $V$ and $-V^*$ have vanishing slow coefficients, and they satisfy the bounds
\[
U_{b_-}\le V\le U_{b_+},\qquad
-U_{b_+}^*\le -V^*\le -U_{b_-}^*.
\]
Since $-U_b^*\in\M_{-b}$,
Lemma~\ref{lem:fast-difference} applies and yields
\[
D^jT=O_\mu(\rho^{-\beta})\qquad\text{for }0\le j\le2.
\]

Now, let $q=S-H$. By Proposition~\ref{prop:coarse}, we have
$D^jq=O_\mu(\rho^{-\alpha})$ for $0\le j\le2$. Since both $q$ and $E$ are odd in $t$, they are
orthogonal to $h$ after multiplication by the even cutoff function from
Lemma~\ref{lem:localization}. We then apply Lemma~\ref{lem:orthogonal} on successively smaller tubular regions.

Substituting $S=H+q$ into the equation for $S$ and using the relation
$\Delta H=H^3-H+Bh$, we obtain the following identity derived in
Appendix~\ref{app:identities}:
\begin{equation*}
\mathcal L_0q
=Bh+\left(A-\frac{n-2}{\rho}\right)q_\rho
+Bq_t-3Hq^2-q^3-3ST^2.
\tag{\ref*{apeq:odd-profile-equation}}
\end{equation*}
Since $H$ is odd, it follows that $q=S-H$ is exactly the odd part of $V-H$.

The term $Bh$ is $O_\mu(\rho^{-2})$, because
$B=O(|t|\rho^{-2})$ and $h$ decays exponentially.
Using the relation $A-(n-2)/\rho=O(t^2\rho^{-3})$ and slightly decreasing $\mu$
to absorb the polynomial growth in $|t|$, we also obtain
\[
\left(A-\frac{n-2}{\rho}\right)q_\rho
=O_\mu(\rho^{-\alpha-3}),\qquad
Bq_t=O_\mu(\rho^{-\alpha-2}).
\]
The terms involving $q^2$ and $q^3$ are $O_\mu(\rho^{-2\alpha})$,
while $ST^2=O_\mu(\rho^{-2\beta})$. Since $\beta>\alpha>1$,
the entire right-hand side is $O_\mu(\rho^{-2})$. Differentiating the right-hand side once and applying the derivative
bounds for $q$ and $T$ yields
$|\nabla(\mathcal L_0q)|=O_\mu(\rho^{-2})$.
Consequently, $\mathcal L_0q=O_\mu(\rho^{-2})$ in $C^{0,\gamma}$
on unit balls, as required for the application of Lemma~\ref{lem:orthogonal}.

We apply \eqref{apeq:orthogonal-estimate} with $k=2$ to $\theta q$, which is odd in $t$ and therefore orthogonal to $h$. The error introduced by the cutoff function is exponentially small by \eqref{apeq:cutoff-smallness}. Consequently, on
the smaller tubular region where $\theta=1$, we obtain
$D^jq=O_\mu(\rho^{-2})$ for $0\le j\le2$.

To estimate $E=q-(n-2)\rho^{-2}\eta$, we use the relation $L_H\eta=-th$ to obtain
\[
\mathcal L_0(\rho^{-2}\eta)
=-\rho^{-2}th+2(n-5)\rho^{-4}\eta.
\]
It follows that
\[
	\begin{aligned}
	\mathcal L_0E
	={}&-\frac{(n-2)t^3}{\rho^2(\rho^2-t^2)}h
	-\frac{2(n-2)(n-5)}{\rho^4}\eta\\
	&+\left(A-\frac{n-2}{\rho}\right)q_\rho+Bq_t
	-3Hq^2-q^3-3ST^2,
	\end{aligned}
\]
where we have also used the cancellation
\[
\left(B+\frac{(n-2)t}{\rho^2}\right)h
=-\frac{(n-2)t^3}{\rho^2(\rho^2-t^2)}h
=O_\mu(\rho^{-4}).
\]
The improved bounds on $q$, the exponential decay of $\eta$ and $h$,
and the bounds on $T$ (noting $\beta\ge3$) establish the same decay order for each
term along with its first derivatives. Hence, $\mathcal L_0E=O_\mu(\rho^{-4})$
in $C^{0,\gamma}$ on unit balls. Finally, applying \eqref{apeq:orthogonal-estimate}
to $\theta E$ yields the desired bounds on $\mathcal T_{\kappa,R}$.
\end{proof}

\section{Ordering property}\label{sec:ordering}
In this section, we show that the minimizers $U_a$ are ordered with respect to the magnitude of $a$, and we prove that such global minimizers possess $O(m)\times O(m)$ symmetry. Finally, we construct the one-sided limits of these minimizers as the parameter of the slow mode tends to zero.

\begin{lemma}\label{lem:comparison}
Let $u,v:\R^n\to(-1,1)$ be global minimizers. If
$(u-v)^+\in L^2(\R^n)$ and $u(x_0) \le v(x_0)$ for some $x_0 \in \mathbb{R}^n$, then either $u\equiv v$ or $u<v$ throughout $\R^n$.

Moreover, if $u=v$ on the boundary of a half-space $\mathcal H$
and $(u-v)^+\in L^2(\mathcal H)$, then $u\le v$ in $\mathcal H$.
\end{lemma}

\begin{proof}
The argument here is inspired by \cite[Proposition~7]{LWW}.
Throughout this proof, all integrals are taken over $\R^n$ unless specified otherwise.
Let $p=(u-v)^+$ and $Q=u^2+uv+v^2-1$.
Subtracting the equations satisfied by $u$ and $v$ yields
$(-\Delta+Q)(u-v)=0$.
Choose a cutoff function $\varphi=\varphi_R\in C_c^\infty(B_{2R})$ satisfying
$0\le\varphi\le1$, $\varphi=1$ on $B_R$, and
$|\nabla\varphi|\le C/R$.
Since $p\in H^1_{\mathrm{loc}}(\mathbb{R}^n)$ and
$\nabla p=\mathbf1_{\{u>v\}}\nabla(u-v)$ almost everywhere,
testing with $\varphi(1-\varphi)p$ gives
\[
\int\varphi(1-\varphi)(|\nabla p|^2+Qp^2)
+\int(1-2\varphi)p\nabla p\cdot\nabla\varphi=0.
\]

The competitors $u_R=u-\varphi p$ and $v_R=v+\varphi p$
are bounded, compactly supported $H^1$ perturbations of $u$ and $v$, so the corresponding energy differences
\[
D_u(R)=\E(u_R;B_{2R})-\E(u;B_{2R}),\qquad
D_v(R)=\E(v_R;B_{2R})-\E(v;B_{2R})
\]
are non-negative. Expanding the gradient terms and combining the pointwise identity
\[
	W(u-\varphi p)+W(v+\varphi p)-W(u)-W(v)
	=-\varphi(1-\varphi)Qp^2 + \tfrac12\varphi^2(1-\varphi)^2p^4
\]
with the relation obtained by testing the equation yields
\begin{equation}
D_u(R)+D_v(R)
=\int p^2|\nabla\varphi|^2
+\tfrac12\int\varphi^2(1-\varphi)^2p^4.
\end{equation}
Both integrands on the right-hand side are supported in $B_{2R}\setminus B_R$.
Since $0\le p\le2$ and $p\in L^2(\R^n)$, letting $R\to\infty$ gives
\begin{equation}\label{eq:energy-difference}
0\le D_u(R)+D_v(R)
\le C(1+R^{-2})\int_{\R^n\setminus B_R}p^2
\to 0.
\end{equation}

Now we construct new global minimizers using $u$ and $v$. Let $w=u-p=\min\{u,v\}$.
Fix $\psi\in C_c^\infty(\Omega)$, where $\Omega$ is bounded,
and take $R$ large enough so that $\overline\Omega\subset B_R$.
Then $u_R=w$ in $\Omega$, and the minimality of $u$ gives
\[
0\le\E(u_R+\psi;B_{2R})-\E(u;B_{2R}) =D_u(R)+\E(w+\psi;\Omega)-\E(w;\Omega).
\]
Letting $R\to\infty$ and applying \eqref{eq:energy-difference} prove that $w$ is a global minimizer.
The same argument shows that $v+p=\max\{u,v\}$ is a global
minimizer. Hence both are smooth solutions.
At any point, at least one of $u-w$ and $v-w$ vanishes, and thus the strong maximum principle gives $w\equiv u$ or $w\equiv v$. This proves the assertion.

For the half-space assertion, set $p=(u-v)^+$ in $\mathcal H$
and extend it by zero to $\R^n\setminus\mathcal H$. The boundary condition ensures that $p$ has zero trace on $\partial\mathcal H$,
so $p\in H^1_{\mathrm{loc}}(\R^n)\cap L^2(\R^n)$.
The same testing argument and exchange identity remain valid, since
$\nabla p=\mathbf1_{\mathcal H\cap\{u>v\}}\nabla(u-v)$
almost everywhere. As before, it follows that
$u-p$ is a global minimizer. Since $u-p\le u$ everywhere and $u-p=u$ on the complementary half-space, the strong maximum
principle yields $u-p\equiv u$, proving that $u\le v$ in $\mathcal H$.
\end{proof}

Using the moving plane method, we prove the following lemma.
\begin{lemma}\label{lem:rotational-symmetry}
Every $u\in\M_a$ with $a\in\R$ is $O(m)\times O(m)$-invariant.
\end{lemma}

\begin{proof}
Fix $e\in\Sph^{m-1}$ and $\lambda>0$. For
$\mathcal H_\lambda=\{x'\cdot e>\lambda\}$, set
$d=x'\cdot e-\lambda$, $x^\lambda=x-2d(e,0)$, and
$v(x)=u(x^\lambda)$.
Combining the identity
\[
t(x)-t(x^\lambda)
=\frac{4\lambda d}
{\sqrt2\bigl(|x'|+|(x^\lambda)'|\bigr)}
\in[0,\sqrt2\lambda]
\]
with Proposition~\ref{prop:coarse} applied to $q=u-H(t)$
yields, on each fixed strip $|t|\le t_0$ and for sufficiently large $\rho$,
\[
H(t(x))-H(t(x^\lambda))\ge c_{\lambda,t_0}\frac d\rho,
\qquad
|q(x)-q(x^\lambda)|
\le C_{\lambda,t_0}\rho^{-\alpha}\min\{d,1\}.
\]
Here, the second estimate follows directly from the pointwise bound on $q$ when $d\ge1$, and by integrating the gradient estimate along the segment joining $x$ and $x^\lambda$ when $d<1$. Since $\alpha>1$, it follows that $u>v$ on this strip outside a sufficiently large ball.

Choose $t_0$ sufficiently large so that Proposition~\ref{prop:away}
ensures $Q=u^2+uv+v^2-1\ge c_0>0$ for $|t|>t_0$.
Then $p=(v-u)^+$ vanishes in the strip $|t|\le t_0$ outside a large ball
and satisfies $(-\Delta+c_0)p\le0$ on its exterior positivity set.
Since $p$ vanishes on the boundary of this region outside the inner sphere,
comparing $p$ with $Ce^{-c|x|}$ via Lemma~\ref{lem:basic-mp}
yields $p\in L^2(\mathcal H_\lambda)$.
Finally, since $u=v$ on $\partial\mathcal H_\lambda$, applying the half-space
assertion of Lemma~\ref{lem:comparison} to $v$ and $u$
gives $u\ge v$ in $\mathcal H_\lambda$.

Letting $\lambda\downarrow0$ and replacing $e$ with $-e$ yields
symmetry across the hyperplane $\{x'\cdot e=0\}$. Since such reflections generate
$O(m)$, it follows that $u$ is radially symmetric in $x'$.
Finally, applying the same argument to $-u(x'',x')\in\M_{-a}$ gives
radial symmetry in $x''$.
\end{proof}

\begin{proposition}
\label{prop:ordering}
If $a<b$, $u\in\M_a$, and $v\in\M_b$, then
$u<v$ in $\R^n$.
\end{proposition}

\begin{proof}
By Lemma~\ref{lem:rotational-symmetry}, both $u$ and $v$
are $O(m)\times O(m)$-invariant. Applying Lemma~\ref{lem:coarse-expansion}
to each solution and subtracting the resulting expansions yields
\[
v-u=c_m(b-a)\rho^{-\alpha}h
+O_\mu(\rho^{-\alpha-\epsilon_*}),
\]
where the leading curvature terms cancel.

Choose $t_0$ sufficiently large that
Proposition~\ref{prop:away} gives
\[
Q:=u^2+uv+v^2-1\ge c_0>0
\qquad\text{for }|t|\ge t_0.
\]
Since $b>a$ and $\min_{|t|\le t_0}h(t)>0$, the preceding
expansion gives $u<v$ on $\{|t|\le t_0\}$ outside a
sufficiently large ball $B_{R_0}$.

Set $p=(u-v)^+$ and
\[
D=\{u>v\}\setminus\overline B_{R_0}.
\]
Then $D\subset\{|t|>t_0\}$, and Kato's inequality for $p=(u-v)^+$ gives
\[
(-\Delta+c_0)p=(c_0-Q)p\le0\quad\text{in }D,
\qquad
p=0\quad\text{on }\partial D\setminus\partial B_{R_0}.
\]
For $0<\mu<\sqrt{c_0}$, the barrier function
$\Phi(x)=2e^{-\mu(|x|-R_0)}$ satisfies
\[
(-\Delta+c_0)\Phi
=\left(c_0-\mu^2+\frac{(n-1)\mu}{|x|}\right)\Phi>0.
\]
Since $p\le2$ on the inner boundary $\partial B_{R_0}$,
applying Lemma~\ref{lem:basic-mp} to $p-\Phi$ yields
\[
p(x)\le2e^{-\mu(|x|-R_0)}
\qquad\text{for }|x|>R_0,
\]
and hence $p\in L^2(\R^n)$.
Finally, Lemma~\ref{lem:comparison} implies that $u<v$ throughout $\R^n$.
\end{proof}

By Proposition~\ref{prop:ordering}, it follows immediately that the limits
\[
U_0^-:=\limsup_{b \uparrow 0}U_b,\qquad
U_0^+:=\liminf_{b \downarrow 0}U_b
\]
are well-defined, independently of the choice of $U_b \in \M_b$.

\begin{lemma}
\label{lem:envelope}
The functions $U_0^-$ and $U_0^+$ belong to $\M_0$ and satisfy
\[
U_0^-=\lim_{b\uparrow0}U_b,\qquad
U_0^+=\lim_{b\downarrow0}U_b
\quad\text{in }C^\infty_{\mathrm{loc}}(\R^n).
\]
Moreover,
\begin{equation}\label{eq:envelope-order}
U_0^-\le V\le U_0^+\qquad\text{for }V\in\M_0,
\end{equation}
and consequently,
\[
U_0^-=\inf_{V\in\M_0}V,\qquad
U_0^+=\sup_{V\in\M_0}V,
\]
where the infimum and supremum are taken pointwise.
\end{lemma}

\begin{proof}
Uniform interior estimates and standard elliptic compactness yield
convergence in $C^\infty_{\mathrm{loc}}(\R^n)$, ensuring that both limits are global minimizers.

For every $0<\delta<\delta_0$, the same ordering property yields
\[
U_{-\delta}\le U_0^-\le U_0^+\le U_\delta.
\]
Consequently, Lemma~\ref{lem:coarse-expansion} applies to both $U_0^-$ and $U_0^+$.
Let $A_\pm$ denote their corresponding slow coefficients. Evaluating the preceding
inequalities at $t=0$, multiplying by $\rho^\alpha/h(0)$,
and sending $\rho\to\infty$, we obtain
\[
-c_m\delta\le A_\pm\le c_m\delta.
\]
Letting $\delta\downarrow0$ then gives $A_\pm=0$.

We next verify that both limits belong to $\M_0$.
Fix $U\in\{U_0^-,U_0^+\}$. By the trapping inequalities and the
decay estimates away from the cone for the reference minimizers, we may choose
$t_0>0$ large enough so that
\[
U(\rho,t)>0\quad\text{for }t\ge t_0,\qquad
U(\rho,t)<0\quad\text{for }t\le-t_0.
\]
Consequently, outside a compact set, the nodal set of $U$ is contained in the fixed strip $|t|<t_0$.

Since the slow coefficient of $U$ is zero,
Lemma~\ref{lem:coarse-expansion} yields
\[
U_t(\rho,t)
=h(t)+(n-2)\rho^{-2}\eta'(t)+o(1)
\]
uniformly on $|t|\le t_0$. For sufficiently large $\rho$, this strip is contained in the tubular region where
the profile estimates hold, so that
\[
U_t(\rho,t)\ge c_0\qquad\text{for }|t|\le t_0.
\]
Combined with $U(\rho,-t_0)<0<U(\rho,t_0)$, this implies that there exists a unique
$f_U(\rho)\in(-t_0,t_0)$ satisfying
\[
U(\rho,f_U(\rho))=0.
\]
By the implicit function theorem, $f_U$ is smooth, and monotonicity gives $\sgn U(\rho,t)=\sgn(t-f_U(\rho))$.

It remains to verify that the parameter of the slow mode of $U$ is zero.
Since $H(0)=\eta(0)=0$ and the slow coefficient vanishes,
Lemma~\ref{lem:coarse-expansion} yields
\[
U(\rho,0)=o(\rho^{-\alpha}).
\]
Combining the mean value theorem with the lower bound on $U_t$ gives $f_U(\rho)=o(\rho^{-\alpha})$.
Since $U$ is $O(m)\times O(m)$-invariant, this estimate is uniform in the angular variables.
All conditions in Definition~\ref{def:moduli} are
therefore satisfied with parameter of the slow mode zero, and hence $U_0^\pm\in\M_0$.

\end{proof}

We conclude this section by recording several properties that will be used in the proof of the pointwise estimates.
\begin{proposition}\label{prop:reflection}
The reflection $(\mathcal RV)(x',x'')=-V(x'',x')$ is an
order-reversing bijection from $\M_a$ onto $\M_{-a}$ satisfying
\[
\mathcal R U_0^\pm=U_0^\mp.
\]

For some $R_0\ge1$ and $\mu\in(0,\sqrt2)$, the two solutions admit
the expansions
\begin{equation}\label{eq:zero-profile}
U_0^\pm=H(t)+(n-2)\rho^{-2}\eta(t)+E^\pm
\qquad\text{in }\mathcal T_{\kappa,R_0},
\end{equation}
where the remainders satisfy
\[
|D^j E^\pm|\le C\rho^{-\min\{4,\beta\}}e^{-\mu|t|}
\qquad (j=0,1,2).
\]
Moreover,
\begin{equation}\label{eq:zero-error}
\begin{aligned}
U_0^\pm-H(t) &= O(R^{-2})
&&\text{as }R\to\infty,\\
D^j\bigl(U_0^\pm+(U_0^\pm)^*\bigr) &= O_\mu(\rho^{-\beta})
&&\text{in }\mathcal T_{\kappa,R_0},\quad j=0,1,2.
\end{aligned}
\end{equation}

Finally, there exist $t_0,\nu>0$ such that, for $|t|\ge t_0$ and $j=0,1,2$,
\[
\bigl|D^j(U_0^\pm-\sgn t)\bigr|
+\bigl|D^j\bigl(U_0^\pm+(U_0^\pm)^*\bigr)\bigr|
\le Ce^{-\nu|t|}.
\]
In particular, on $\{\rho>R_0,\ |t|\ge\kappa\rho\}$, the right-hand side is bounded by
$Ce^{-cR}$, where $c=\nu\kappa/\sqrt{1+\kappa^2}$.
\end{proposition}

\begin{proof}
Since the potential is even, minimality is preserved under orthogonal transformations and sign changes.
Moreover, the reflected nodal set is a graph given by $-f_V(\rho,\theta,\omega)$, which has the slow mode parameter $-a$ and the prescribed phase orientation.
Thus, $\mathcal R$ maps $\M_a$ into $\M_{-a}$.
Since $\mathcal R^2=\mathrm{Id}$ and $\mathcal R$ is order-reversing by definition, it follows that $\mathcal R$ is an order-reversing bijection from $\M_a$ onto $\M_{-a}$.

By Lemma~\ref{lem:envelope}, $U_0^-$ and $U_0^+$ are the
pointwise infimum and supremum of $\M_0$, respectively.
Since $\mathcal R$ is an order-reversing bijection of $\M_0$ onto itself, we have
\begin{align*}
\mathcal RU_0^+
&=\inf_{V\in\M_0}\mathcal RV
=\inf_{V\in\M_0}V=U_0^-,\\
\mathcal RU_0^-
&=\sup_{V\in\M_0}\mathcal RV
=\sup_{V\in\M_0}V=U_0^+.
\end{align*}

Fix $U\in\{U_0^-,U_0^+\}$. By Lemma~\ref{lem:envelope},
$U \in \M_0$ is trapped between
$U_{-\delta}$ and $U_\delta$ for any fixed small $\delta>0$.
Thus, setting $S=(U-U^*)/2$ and $T=(U+U^*)/2$,
Corollary~\ref{cor:sharp-profile} yields
\[
D^j\bigl(S-H-(n-2)\rho^{-2}\eta\bigr)=O_\mu(\rho^{-4}),
\qquad D^jT=O_\mu(\rho^{-\beta}),
\qquad 0\le j\le2.
\]
Since $U=S+T$ and $U+U^*=2T$, combining these estimates with $\beta\ge3$,
the exponential decay of $\eta$ from Lemma~\ref{lem:normal},
and the comparability of $\rho$ and $R$ in the tubular region establishes
\eqref{eq:zero-profile} and \eqref{eq:zero-error}.

For the estimates away from the cone, Proposition~\ref{prop:away} yields
\[
|U-\sgn t|\le Ce^{-\nu|t|}
\qquad\text{for }|t|\ge t_0.
\]
On any ball $B_2(x)$ sufficiently far from the cone,
$\sigma=\sgn t(x)$ is constant and $U-\sigma$ satisfies
\[
(-\Delta+U^2+\sigma U)(U-\sigma)=0.
\]
The standard interior elliptic estimates yield
\[
\sum_{j=0}^2
\left(
|D^j(U-\sgn t)|+|D^j(U+U^*)|
\right)
\le Ce^{-\nu|t|}
\qquad\text{for }|t|\ge t_0.
\]
\end{proof}

\section{Maximum principles and pointwise estimates}
\label{sec:signs}

Throughout this section, all functions and test functions are assumed to be
$O(m)\times O(m)$-invariant. We use the standard Euclidean differential
operators and volume measure $dx$, and work on the domain
\[
\Omega=\{x=(x',x'')\in\R^m\times\R^m:|x'|>|x''|>0\},
\]
with $r=|x'|$, $s=|x''|$, and $R=|x|$.
Several identities used below are derived in Appendix~\ref{app:identities}.

Our arguments are closely related to
the derivative estimates of
Cabr\'e--Terra \cite[Theorem~1.7]{CabreTerra2010}
and Cabr\'e \cite[Proposition~1.5 and Section~4]{Cabre2012}
for saddle-shaped solutions.
Because odd symmetry across the cone is not assumed in our setting,
we adapt their approach under the present stability and
asymptotic hypotheses, and then establish the coupled
derivative inequalities for the odd and even parts.

We first establish a maximum principle for non-negative elliptic
operators on unbounded domains under suitable decay assumptions
at infinity.

\begin{lemma}
\label{lem:stable-mp}
Let $m\ge4$ and set
\[
K=(m-1)\left(\frac1{r^2}+\frac1{s^2}\right).
\]
Let $q$ be a bounded scalar function satisfying
\[
\int_\Omega\bigl(|\nabla\varphi|^2+q\varphi^2\bigr)\,dx
\ge0 \qquad\text{for all }\varphi\in C_c^\infty(\Omega).
\]
Suppose that $z\in C^2(\Omega)\cap C^{0,1}_{\mathrm{loc}}
(\overline\Omega)$ satisfies
\[
(-\Delta+q+K)z\ge0\quad\text{in }\Omega,
\qquad z=0\quad\text{on }\{r=s>0\},
\]
along with the decay condition
\[
|z|=O(R^{1-\beta})\qquad \text{as }R\to\infty.
\]
Then $z\ge0$ in $\Omega$. Moreover, if the differential inequality is strict, then $z>0$ in $\Omega$.
\end{lemma}

\begin{proof}
To ensure that the cutoff error at infinity vanishes under the assumed decay of $z$, we introduce the auxiliary weight
$\psi=(1+R^2)^{-1}$ and choose a
smooth radial cutoff function $0\le\theta_N\le1$ satisfying
\[
\theta_N=1\quad\text{on }\{R\le N\},\qquad
\theta_N=0\quad\text{on }\{R\ge2N\},\qquad
|\nabla\theta_N|\le C/N.
\]

We first justify testing the equation near the axis $\{s=0\}$.
Let $\chi_\epsilon(s)$ be a smooth cutoff function equal to $0$ for
$s\le\epsilon$ and $1$ for $s\ge2\epsilon$, with
$|\nabla\chi_\epsilon|\le C/\epsilon$.
Fix $N$ and set $\varphi=\psi\theta_N$.
Testing the differential inequality against
$\varphi^2\chi_\epsilon^2z_-$ and using
\[
\nabla z_-\cdot\nabla(\varphi^2\chi_\epsilon^2z_-)
= |\nabla(\varphi\chi_\epsilon z_-)|^2
-|\nabla(\varphi\chi_\epsilon)|^2z_-^2,
\]
we obtain
\[
	\begin{aligned}
	&\int_\Omega
	\left(|\nabla(\varphi\chi_\epsilon z_-)|^2
	+q\varphi^2\chi_\epsilon^2z_-^2	\right)dx
	+\int_\Omega K\varphi^2\chi_\epsilon^2z_-^2\,dx \le
	\int_\Omega|\nabla(\varphi\chi_\epsilon)|^2z_-^2\,dx.
	\end{aligned}
\]
By our positivity assumption applied to $\varphi\chi_\epsilon z_-$, the first integral is non-negative, and hence
\[
\int_\Omega K\varphi^2\chi_\epsilon^2z_-^2\,dx
\le \int_\Omega|\nabla(\varphi\chi_\epsilon)|^2z_-^2\,dx.
\]

We next let $\epsilon\downarrow0$ with $N$ fixed.
Since $s=|x''|$ and $x''\in\R^m$, we have
\[
\bigl|B_{2N}\cap\{s<2\epsilon\}\bigr|
\le C_N\epsilon^m.
\]
Combined with the local boundedness of $z_-$ and the gradient bound on
$\chi_\epsilon$, this yields
\[
\int_\Omega\varphi^2|\nabla\chi_\epsilon|^2z_-^2\,dx
\le C_N\epsilon^{m-2}\to 0,
\]
as $m>2$. Moreover, by the Cauchy--Schwarz inequality,
\[
	\begin{aligned}
	&\left| 2\int_\Omega \varphi\chi_\epsilon z_-^2
	\nabla\varphi\cdot\nabla\chi_\epsilon\,dx \right|\\
	&\qquad\le
	2\left(\int_\Omega\chi_\epsilon^2|\nabla\varphi|^2z_-^2\,dx\right)^{1/2}
	\left(\int_\Omega\varphi^2|\nabla\chi_\epsilon|^2z_-^2\,dx\right)^{1/2}
	\to 0.
	\end{aligned}
\]
Expanding $|\nabla(\varphi\chi_\epsilon)|^2$ and applying the
dominated convergence theorem to the remaining term, we conclude that
\[
\int_\Omega|\nabla(\varphi\chi_\epsilon)|^2z_-^2\,dx
\to
\int_\Omega|\nabla\varphi|^2z_-^2\,dx.
\]

Since $K\ge0$ and $\chi_\epsilon \to 1$ pointwise in $\Omega$,
applying Fatou's lemma to the left-hand side yields
\[
\int_\Omega K\psi^2\theta_N^2z_-^2\,dx
\le
\int_\Omega|\nabla(\psi\theta_N)|^2z_-^2\,dx.
\]
Using $|\nabla\psi|^2\le4R^{-2}\psi^2$ together with Young's inequality, we find
\[
	\begin{aligned}
	|\nabla(\psi\theta_N)|^2
	&\le \frac54\theta_N^2|\nabla\psi|^2 
	+10\psi^2|\nabla\theta_N|^2\\
	&\le 5R^{-2}\psi^2\theta_N^2
	+10\psi^2|\nabla\theta_N|^2.
	\end{aligned}
\]

Combining $K\ge4(m-1)R^{-2}\ge 6R^{-2}$ with the decay assumption on $z$
then yields
\begin{equation}\label{eq:absorption}
	\begin{aligned}
	\int_\Omega R^{-2}\psi^2\theta_N^2z_-^2\,dx
	&\le 10\int_\Omega\psi^2|\nabla\theta_N|^2z_-^2\,dx\\
	&\le CN^{n-4-2\beta}.
	\end{aligned}
\end{equation}
Here $n-4-2\beta\le 6-n<0$, since $\beta\ge n-5$ and $n\ge8$.
Letting $N\to\infty$ shows that $z_-=0$, and hence $z\ge0$ in $\Omega$.
Finally, if the differential inequality is strict, any interior zero of $z$
would be a local minimum at which $(-\Delta+q+K)z\le0$,
a contradiction.
\end{proof}

For a stable solution $U$, the second variation
\[
Q_U[\xi]=\int_{\R^n}
\bigl(|\nabla\xi|^2+(3U^2-1)\xi^2\bigr)\,dx
\]
satisfies $Q_U[\xi]\ge0$ for all compactly supported $\xi\in H^1(\R^n)$.

\begin{lemma}
\label{lem:coordinate-signs}
Every $U\in\M_0$ satisfies
\[
U_r>0\quad\text{for }r>0,\qquad
U_s<0\quad\text{for }s>0.
\]
\end{lemma}

\begin{proof}
By Lemma~\ref{lem:rotational-symmetry}, $U$ is
$O(m)\times O(m)$-invariant. Moreover, global minimality implies stability,
and Propositions~\ref{prop:coarse} and~\ref{prop:away}
yield \eqref{eq:coarse} and \eqref{eq:away}, respectively.

By Lemma~\ref{lem:derivatives}, the derivatives $U_r$ and $U_s$ satisfy
\begin{equation*}
\left(\mathcal L_U+\frac{m-1}{r^2}\right)U_r=0,\qquad
\left(\mathcal L_U+\frac{m-1}{s^2}\right)U_s=0,
\tag{\ref*{apeq:radial-linearizations}}
\end{equation*}
where $\mathcal L_U=-\Delta+3U^2-1$.

On every fixed normal strip $\{|t|\le T\}$, we have
$U_r\to h(t)/\sqrt2$ and $-U_s\to h(t)/\sqrt2$ uniformly
as $\rho\to\infty$.
We may choose $t_0>0$ large enough so that
$3U^2-1\ge\kappa_0>0$ on $\{|t|\ge t_0\}$ outside a large ball.
Since $h$ has a positive minimum on $[-t_0-1,t_0+1]$,
we can then choose $R_0$ sufficiently large so that
$U_r>0$ and $-U_s>0$ in
$\{R>R_0,\ |t|\le t_0+1\}$.

Setting $p=(U_r)_-$, it follows that $p=0$ in $\{R>R_0,\ |t|\le t_0+1\}$,
whereas $3U^2-1\ge\kappa_0$ whenever $R>R_0$ and $|t|>t_0$.
By Kato's inequality, $p$ satisfies
$(-\Delta+\kappa_0)p\le0$ weakly in $\{R>R_0\}$,
and the standard maximum principle then yields
$p\le Ce^{-c|x|}$.

Near the axis $\{r=0\}$, we have $U_r=O(r)$.
Using the cutoff functions from the proof of
Lemma~\ref{lem:stable-mp}, we test the differentiated equation for $U_r$
against $\theta_N^2\chi_\epsilon(r)^2p$ and invoke stability.
Removing the cutoff near $r=0$ by the same argument (with $r$ in place of $s$),
and then letting $N\to\infty$ via the exponential decay of $p$, we obtain
\[
Q_U[p]+(m-1)\int_{\R^n}r^{-2}p^2\,dx\le0.
\]
Since both terms on the left-hand side are non-negative, we conclude that $p=0$, and the strong maximum
principle yields $U_r>0$ for $r>0$. The argument for $-U_s$ is identical.
\end{proof}

We now prove the key estimate.

\begin{proposition}[Pointwise estimates]\label{prop:pointwise-estimates}
Let $U = U_0^+$. Set $x^*=(x'',x')$,
$U^*(x)=U(x^*)$, and
\[
S=\frac{U-U^*}{2},\qquad T=\frac{U+U^*}{2}.
\]
Fix $0<\kappa<1/4$. There exist $R_0,C,\mu,c>0$
such that the following hold.

\begin{enumerate}[label=\textup{(\roman*)}]
\item \textup{Estimates in the tubular region:}
For $\rho>R_0$ and $|t|<\kappa\rho$,
\begin{equation}\label{eq:stress-derivative-bounds}
	\begin{aligned}
	|S_\rho|+|S_{\rho t}|
	&\le C\rho^{-3}e^{-\mu|t|},\\
	|T|+|\nabla T|
	&\le C\rho^{-\beta}e^{-\mu|t|}.
	\end{aligned}
\end{equation}

\item \textup{Estimates away from the cone:}
For $\rho>R_0$ and $|t|\ge\kappa\rho$,
\begin{equation}\label{eq:stress-away-bounds}
|S_\rho|+|S_{\rho t}|+|T|+|\nabla T|
\le Ce^{-cR}.
\end{equation}

\item \textup{Signs of derivatives:}
If $T>0$, then for $0<t<\rho$,
\begin{equation}\label{eq:all-signs}
S>0,\quad T>0,\quad S_\rho>0,\quad S_t>0,
\quad T_\rho<0,\quad T_t<0.
\end{equation}
\end{enumerate}
\end{proposition}

\begin{proof}
By Lemma~\ref{lem:envelope}, $U$ is an
$O(m)\times O(m)$-invariant global minimizer and is therefore stable.
Moreover, Proposition~\ref{prop:reflection} yields the asymptotic estimates
\eqref{eq:zero-profile}--\eqref{eq:zero-error}.
In particular, the estimate for $|T|+|\nabla T|$ in \textup{(i)} follows immediately from
\eqref{eq:zero-error}.

For $S$, Corollary~\ref{cor:sharp-profile} yields
$S=H+(n-2)\rho^{-2}\eta+E$, with
$|E_\rho|+|E_{\rho t}|=O_\mu(\rho^{-4})$. Since $H$ depends
only on $t$, we have
\[
S_\rho=-2(n-2)\rho^{-3}\eta+E_\rho,\qquad
S_{\rho t}=-2(n-2)\rho^{-3}\eta'+E_{\rho t}.
\]
Combining these identities with the exponential decay of $\eta$ and $\eta'$ establishes the bound
for $S$ in \textup{(i)}.

To prove \textup{(ii)}, we combine exponential convergence to the constant
phases with standard interior elliptic estimates to obtain
\[
|S_\rho|+|S_{\rho t}|+|T|+|\nabla T|
\le Ce^{-\nu|t|}
\]
for sufficiently large $|t|$.
Since $R^2=\rho^2+t^2$ and $|t|\le\rho$, the condition
$|t|\ge\kappa\rho$ implies $|t|\ge\kappa R/\sqrt2$,
which immediately yields \eqref{eq:stress-away-bounds}.

We now prove \textup{(iii)} under the assumption $T>0$.

\medskip
\textbf{STEP 1}: We have $\mathcal K S, -\mathcal K T>0$ in $\Omega$.

\begin{proof}

Lemma~\ref{lem:coordinate-signs} implies that
$S_r>0$ and $S_s<0$ in $\Omega$, so $\mathcal K S=sS_r-rS_s>0$.
Since $S=0$ on the cone, we also obtain $S>0$ and $S_t>0$
in $\Omega$.
Recall that the odd and even parts satisfy
\begin{equation*}
-\Delta S=S(1-S^2-3T^2),\qquad
-\Delta T=T(1-T^2-3S^2),
\tag{\ref*{apeq:ST-equations}}
\end{equation*}
and the differentiated identities used below are derived in
Lemma~\ref{lem:derivatives}.
Define
\[
\mathcal A=-\Delta+3(S^2+T^2)-1,
\qquad \mathcal K=s\partial_r-r\partial_s.
\]
Since
\[
3(S^2+T^2)-1
=\frac12\bigl((3U^2-1)+(3(U^*)^2-1)\bigr),
\]
the operator $\mathcal A$ is non-negative by the
stability of $U$ and $U^*$.

Applying $\mathcal K$ to the equation for $T$ and using
the commutator identity \eqref{apeq:commutators}, we obtain
\begin{equation*}
\left(\mathcal A+(m-1)(r^{-2}+s^{-2})\right)(-\mathcal K T)
=6ST\,\mathcal K S.
\tag{\ref*{apeq:angular-T-equation}}
\end{equation*}
The right-hand side is strictly positive in $\Omega$, since $S,T>0$ and
$\mathcal K S>0$.
Moreover, the evenness of $T$ implies $\mathcal K T=0$
on $\{r=s>0\}$, while \textup{(i)} and \textup{(ii)} yield
$\mathcal K T=O(R^{1-\beta})$ as $R\to\infty$.
Applying Lemma~\ref{lem:stable-mp} therefore gives
\begin{equation}\label{eq:angular-sign}
-\mathcal K T>0\qquad\text{in }\Omega.
\end{equation}
\end{proof}

The radial derivatives $S_\rho$ and $T_\rho$ satisfy the coupled system
\begin{equation*}
	\begin{aligned}
	\left(\mathcal A+\frac{m-1}{rs}\right)S_\rho
	+6ST\,T_\rho
	&=\frac{(m-1)(r-s)}{\sqrt2\,r^2s^2}\,\mathcal K S,\\
	\left(\mathcal A+\frac{m-1}{rs}\right)(-T_\rho)
	-6ST\,S_\rho
	&=-\frac{(m-1)(r-s)}{\sqrt2\,r^2s^2}\,\mathcal K T,
	\end{aligned}
\tag{\ref*{apeq:coupled-radial-equations}}
\end{equation*}
where both right-hand sides are strictly positive in $\Omega$ by Step~1.
Moreover, the oddness of $S$ and evenness of $T$ yield the boundary conditions
\begin{equation}\label{eq:mixed-boundary}
S_\rho=0,\qquad \partial_\nu T_\rho=0
\qquad\text{on }\{r=s>0\}.
\end{equation}

We next verify the non-negativity of the quadratic form
associated with this coupled system.
Given compactly supported, $O(m)\times O(m)$-invariant test functions
$\varphi,\psi$ on $\overline\Omega$ with $\varphi=0$ on the cone $\{r=s\}$,
we define
\[
\xi(x)=
	\begin{cases}
	\psi(x)-\varphi(x),&r>s,\\
	\psi(x^*)+\varphi(x^*),&r<s.
	\end{cases}
\]
Since $x^*=x$ and $\varphi=0$ on the cone, the traces from both sides agree across $\{r=s\}$,
so $\xi$ is an admissible $H^1(\R^n)$ test function.
The stability of $U$ then yields
\begin{equation}\label{eq:system-stability}
	\begin{aligned}
	0\le\frac12 Q_U[\xi]
	=\int_\Omega\Bigl(
	&|\nabla\varphi|^2+|\nabla\psi|^2\\
	&+\bigl(3(S^2+T^2)-1\bigr)(\varphi^2+\psi^2)
	-12ST\varphi\psi\Bigr)\,dx.
	\end{aligned}
\end{equation}

Before testing the coupled system, we first establish the decay
of $(S_\rho)_-$ and $(T_\rho)_+$.

\medskip
\textbf{STEP 2}: We prove the decay estimate $(S_\rho)_-+(T_\rho)_+=O(R^{-\beta})$.

\begin{proof}
First, note that the bound for $(T_\rho)_+$ follows directly from \textup{(i)} and \textup{(ii)}, so it suffices to estimate $(S_\rho)_-$. 

For the correction function $\eta$,
Lemma~\ref{lem:normal} yields
$\eta(t)<0$ for $t>0$ and $\eta'(0)<0$.
Consequently, for each fixed $t_0>0$, there exists $c>0$ such that
\[
-\eta(t)\ge ct\qquad\text{for }0\le t\le t_0.
\]
Since $E$ is odd in $t$, we have $E_\rho(\rho,0)=0$, and integrating
$E_{\rho t}=O_\mu(\rho^{-4})$ from $0$ to $t$ gives
\[
|E_\rho(\rho,t)|\le Ct\rho^{-4}
\qquad\text{for }0\le t\le t_0.
\]
Hence,
\[
S_\rho=-2(n-2)\rho^{-3}\eta+E_\rho
\ge t\rho^{-3}\left(2(n-2)c-\frac{C}{\rho}\right)>0
\]
for $0<t\le t_0$ and all sufficiently large $\rho$.

By exponential convergence to the constant phases, we may
choose $t_0$ sufficiently large so that the potential in
$\mathcal A+(m-1)/(rs)$ is bounded below by a positive
constant for $t\ge t_0$ outside a large ball.
Increasing $R_0$ if necessary, we deduce from Kato's inequality and the first equation in
\eqref{apeq:coupled-radial-equations} that
\[
	\begin{aligned}
	\left(\mathcal A+\frac{m-1}{rs}\right)(S_\rho)_-
	&\le-\mathbf1_{\{S_\rho<0\}}
	\left(\mathcal A+\frac{m-1}{rs}\right)S_\rho\\
	&\le6ST\,\mathbf1_{\{S_\rho<0\}}T_\rho\\
	&\le6ST(T_\rho)_+
	\end{aligned}
\]
in $\{\rho>R_0,\ t_0<t<\kappa\rho\}$.
Note that $(S_\rho)_-=0$ on $\{t=t_0\}$, while \textup{(ii)} controls
its values on $\{t=\kappa\rho\}$. Moreover, \textup{(i)} yields
$
(T_\rho)_+\le C\rho^{-\beta}e^{-\mu t}.
$

Therefore, we have
\[
(-\Delta+c_0)(S_\rho)_- \le C\rho^{-2\beta}e^{-2\mu t}.
\]
Decreasing $\mu$ and $c$ if necessary, a comparison argument with a sufficiently
large multiple of
$
\rho^{-2\beta}e^{-\mu t}+e^{-cR}
$ yields
\[
(S_\rho)_-
\le C\rho^{-2\beta}e^{-\mu t}+Ce^{-cR}.
\]
Combined with \textup{(i)}, \textup{(ii)}, and the bound
$\rho\ge R/\sqrt2$, this establishes
\[
(S_\rho)_-+(T_\rho)_+=O(R^{-\beta}).
\]

\end{proof}

\textbf{STEP 3:} We apply the maximum principle to the coupled system.

\begin{proof}
Set $f=S_\rho$ and $g=-T_\rho$.
We now adapt the cutoff argument from
Lemma~\ref{lem:stable-mp} to the coupled system \eqref{apeq:coupled-radial-equations}. Since $\mathcal K S>0$ and $\mathcal K T<0$ in $\Omega$, both right-hand sides of this system are strictly positive in $\Omega$.
Setting $\psi_N=\frac{\theta_N}{1+R^2}$,
we test the two equations against
$\psi_N^2f_-$ and $\psi_N^2g_-$, respectively,
initially inserting the axis cutoff $\chi_\epsilon$ from Lemma~\ref{lem:stable-mp}.
The boundary terms on the cone $\{r=s>0\}$ vanish since
$f_-=0$ and $\partial_\nu g=0$ there.
Moreover, because $f_-$ and $g_-$ are locally bounded,
we can pass to the limit $\epsilon\downarrow0$ to remove the axis cutoff just as in Lemma~\ref{lem:stable-mp}.

We obtain
\[
	\begin{aligned}
	&\int_\Omega\Bigl[
	|\nabla(\psi_Nf_-)|^2+|\nabla(\psi_Ng_-)|^2\\
	&\qquad+
	\left(3(S^2+T^2)-1+\frac{m-1}{rs}\right)
	\psi_N^2(f_-^2+g_-^2)
	\Bigr]\,dx\\
	&\le
	\int_\Omega|\nabla\psi_N|^2(f_-^2+g_-^2)\,dx
	-6\int_\Omega ST\psi_N^2(gf_-+fg_-)\,dx.
	\end{aligned}
\]

Applying the stability inequality \eqref{eq:system-stability} to $(\varphi,\psi)=(\psi_Nf_-,\psi_Ng_-)$ and using $ST\ge0$ together with the pointwise bound
\[
gf_-+fg_-=-2f_-g_-+g_+f_-+f_+g_- \ge -2f_- g_-,
\]
we obtain
\[
\int_\Omega\frac{m-1}{rs}
\frac{\theta_N^2(f_-^2+g_-^2)}{(1+R^2)^2}\,dx
\le
\int_\Omega
\left|\nabla\left(\frac{\theta_N}{1+R^2}\right)\right|^2
(f_-^2+g_-^2)\,dx.
\]
Since $rs\le R^2/2$ and $m\ge4$, we have
\[
\frac{m-1}{rs}\ge\frac{2(m-1)}{R^2}\ge\frac{6}{R^2}.
\]
Repeating the absorption argument from
\eqref{eq:absorption} then yields
\[
	\begin{aligned}
	\int_\Omega
	\frac{\theta_N^2(f_-^2+g_-^2)}{R^2(1+R^2)^2}\,dx
	&\le 10\int_\Omega
	\frac{|\nabla\theta_N|^2(f_-^2+g_-^2)}{(1+R^2)^2}\,dx\\
	&\le CN^{n-6-2\beta}.
	\end{aligned}
\]

Here $n-6-2\beta\le4-n<0$, since $\beta\ge n-5$.
Letting $N\to\infty$ yields $f_-=g_-=0$,
and hence $S_\rho\ge0$ and $T_\rho\le0$ in $\Omega$.
Since both right-hand sides of
\eqref{apeq:coupled-radial-equations} are strictly positive in $\Omega$,
the strong maximum principle gives the strict inequalities
\[
S_\rho>0,\qquad T_\rho<0
\qquad\text{in }\Omega.
\]
\end{proof}

Finally, $\mathcal K=-t\partial_\rho+\rho\partial_t$, thus
$
\rho T_t=\mathcal K T+tT_\rho<0
$
in $\Omega$  by \eqref{eq:angular-sign}.
This proves $T_t<0$ and completes the proof.
\end{proof}

\section{Proof of the main theorem}
\label{sec:main}

We begin by introducing the relative stress--energy tensor and deriving its conservation law. In what follows, $\omega=(\rho^2-t^2)^{m-1}$ denotes, up to a dimensional constant, the Euclidean volume density in
$(\rho,t)$ coordinates after integrating out the angular variables.

\begin{lemma}[Relative stress--energy tensor]\label{lem:stress}
Let $U$ be an $O(m)\times O(m)$-invariant solution of \eqref{eq:AC},
and set $S=(U-U^*)/2$ and $T=(U+U^*)/2$.
Given the Allen--Cahn energy density
\[
e(u)=\frac12|\nabla u|^2+\frac14(1-u^2)^2,
\]
we define the relative energy density by
\begin{equation}\label{eq:interaction-density}
	e_{\mathrm{rel}}:=\frac{e(-U^*)-e(U)}2=-S_\rho T_\rho-S_tT_t+ST(1-S^2-T^2).
\end{equation}

Then, in the region $\{0<t<\rho\}$, we have
\begin{align}
\partial_\rho\!\left[\omega(2S_\rho T_\rho+e_{\mathrm{rel}})\right]
+\partial_t\!\left[\omega(S_\rho T_t+S_tT_\rho)\right]
&=\omega_\rho e_{\mathrm{rel}},
\label{eq:stress-rho}\\
\partial_\rho\!\left[\omega(S_\rho T_t+S_tT_\rho)\right]
+\partial_t\!\left[\omega(2S_tT_t+e_{\mathrm{rel}})\right]
&=\omega_t e_{\mathrm{rel}}.
\label{eq:stress-t}
\end{align}
\end{lemma}

\begin{proof}
Following \cite[Section~2.1]{Alikakos2011}, we consider the Euclidean stress--energy tensor
\[
\mathsf T[u]=\nabla u\otimes\nabla u-e(u)I_n.
\]
For any solution $u$ of the Allen--Cahn equation \eqref{eq:AC}, a direct computation gives the divergence-free identity
\[
\mathrm{div}\mathsf T[u] =(\Delta u+u-u^3)\nabla u=0.
\]

Since $U=S+T$ and $-U^*=S-T$, expanding $e(-U^*)$ and $e(U)$ yields
\eqref{eq:interaction-density}. The corresponding
relative stress--energy tensor is given by
\[
	\begin{aligned}
	\mathsf T_{\mathrm{rel}}
	&:=\frac{\mathsf T[-U^*]-\mathsf T[U]}2\\
	&=-\nabla S\otimes\nabla T-\nabla T\otimes\nabla S
	-e_{\mathrm{rel}}I_n.
	\end{aligned}
\]
Since both $U$ and $-U^*$ solve the Allen--Cahn equation, it follows that
\[
\mathrm{div}\mathsf T_{\mathrm{rel}}=0.
\]

For $r,s>0$, choose local orthonormal frames
$\{X_k\}_{k=1}^{m-1}$ and $\{Y_k\}_{k=1}^{m-1}$ tangent to
the spheres in the $x'$- and $x''$-variables, respectively.
Together with $\nabla\rho$ and $\nabla t$, these form a local
orthonormal frame on $\R^n$.
By $O(m)\times O(m)$-invariance, $\nabla S$
and $\nabla T$ lie in $\mathrm{span}\{\nabla\rho,\nabla t\}$. Consequently, the
angular components of $\mathsf T_{\mathrm{rel}}$ are given by
\[
\mathsf T_{\mathrm{rel}}(X_k,X_\ell)
=\mathsf T_{\mathrm{rel}}(Y_k,Y_\ell)
=-e_{\mathrm{rel}}\delta_{k\ell},
\]
while all mixed components between the two angular distributions,
as well as between either angular distribution and $\mathrm{span}\{\nabla\rho,\nabla t\}$,
vanish.

Since $\rho=(r+s)/\sqrt2$ and $t=(r-s)/\sqrt2$, their
angular Hessian components are
\[
	\begin{aligned}
	D^2\rho(X_k,X_\ell)&=\frac{\delta_{k\ell}}{\sqrt2\,r},
	&D^2\rho(Y_k,Y_\ell)&=\frac{\delta_{k\ell}}{\sqrt2\,s},\\
	D^2t(X_k,X_\ell)&=\frac{\delta_{k\ell}}{\sqrt2\,r},
	&D^2t(Y_k,Y_\ell)&=-\frac{\delta_{k\ell}}{\sqrt2\,s}.
	\end{aligned}
\]

Since both Hessians vanish on $\mathrm{span}\{\nabla\rho,\nabla t\}$,
only the angular components contribute to the contraction for $i=\rho,t$:
\[
	\begin{aligned}
	\mathsf T_{\mathrm{rel}}:D^2i
	&=-e_{\mathrm{rel}}\sum_{k=1}^{m-1}
	\bigl[D^2i(X_k,X_k)+D^2i(Y_k,Y_k)\bigr]\\
	&=-e_{\mathrm{rel}}\Delta i.
	\end{aligned}
\]
Combining this with $\mathrm{div}\mathsf T_{\mathrm{rel}}=0$,
we obtain
\[
\mathrm{div}(\mathsf T_{\mathrm{rel}}\nabla i)
=\mathsf T_{\mathrm{rel}}:D^2i
=-e_{\mathrm{rel}}\Delta i
=-\frac{\omega_i}{\omega}e_{\mathrm{rel}}.
\]

For $i,j\in\{\rho,t\}$, write
\[
	\begin{aligned}
	(\mathsf T_{\mathrm{rel}})_{ij}
	&:=\mathsf T_{\mathrm{rel}}(\nabla i,\nabla j)\\
	&=-S_iT_j-T_iS_j-\delta_{ij}e_{\mathrm{rel}}.
	\end{aligned}
\]
Since the mixed radial-angular components vanish, we have
\[
\mathsf T_{\mathrm{rel}}\nabla i
=(\mathsf T_{\mathrm{rel}})_{i\rho}\nabla\rho
+(\mathsf T_{\mathrm{rel}})_{it}\nabla t.
\]
Using $\Delta\rho=\omega_\rho/\omega$ and
$\Delta t=\omega_t/\omega$ to compute the divergence of this vector field,
we find
\[
\partial_\rho\!\left[\omega(\mathsf T_{\mathrm{rel}})_{i\rho}\right]
+\partial_t\!\left[\omega(\mathsf T_{\mathrm{rel}})_{it}\right]
=-\omega_i e_{\mathrm{rel}}.
\]
Finally, setting $i=\rho$ and $i=t$, substituting the components above,
and multiplying by $-1$, we obtain
\eqref{eq:stress-rho} and \eqref{eq:stress-t},
respectively.
\end{proof}

\begin{lemma}\label{lem:rigidity}
We have $U_0^+ - U_0^- = 0$ in $\mathbb R^n$.
\end{lemma}

\begin{proof}
We set $U= U_0^+$ and adopt the notation of Lemma~\ref{lem:stress}. Since $U+U^*\ge0$, we have $T \ge 0$.
Suppose toward a contradiction that $T\not\equiv0$. Then the linear equation satisfied by $T$ and the strong
maximum principle yield $T>0$.
Moreover, since $S^2+T^2<1$, Proposition~\ref{prop:pointwise-estimates}\textup{(iii)}
implies that
\[
e_{\mathrm{rel}}>0,\qquad S_\rho T_\rho<0
\quad\text{on }0<t<\rho.
\]

Define
\[
P(\rho)
=-2\int_0^\rho\int_0^{\rho'}
\omega_t(\rho',t')\,
e_{\mathrm{rel}}(\rho',t')\,dt'\,d\rho',
\]
where
\[
\omega_t(\rho',t')
=-(n-2)t'\bigl((\rho')^2-(t')^2\bigr)^{m-2}.
\]
Then $P(0)=0$ and
\begin{equation}\label{eq:Pprime}
P'(\rho)
=-2\int_0^\rho\omega_t e_{\mathrm{rel}}\,dt
=2(n-2)\int_0^\rho
t(\rho^2-t^2)^{m-2}e_{\mathrm{rel}}\,dt>0.
\end{equation}

Integrating the identity
\eqref{eq:stress-t} in $t$ from $0$ to $\rho$ yields
\[
\frac{d}{d\rho}
\left[-2\int_0^\rho
\omega(S_\rho T_t+S_tT_\rho)\,dt\right]
=P'(\rho),
\]
where the boundary terms vanish since
$S=S_\rho=T_t=e_{\mathrm{rel}}=0$ at $t=0$ and
$\omega=0$ at $t=\rho$.
Since the integral on the left-hand side tends to zero as $\rho\downarrow0$, integrating in $\rho$ gives
\[
P(\rho)
=-2\int_0^\rho
\omega(S_\rho T_t+S_tT_\rho)\,dt.
\]

Since $\mathrm{div}\mathsf T_{\mathrm{rel}}=0$ and
$D^2t$ vanishes on the plane spanned by $\nabla S$ and
$\nabla T$, we have
\[
\mathrm{div}(\mathsf T_{\mathrm{rel}}\nabla t)
=\mathsf T_{\mathrm{rel}}:D^2t
=-e_{\mathrm{rel}}\Delta t
=-\frac{\omega_t}{\omega}e_{\mathrm{rel}}.
\]
Consequently, the double integral defining $P$ represents, up to a dimensional constant, twice the volume
integral of this divergence.

Applying the divergence theorem together with the vanishing of the boundary terms
above therefore yields
\[
P(\rho)=2\int_0^\rho
\omega\,\mathsf T_{\mathrm{rel}}(\nabla\rho,\nabla t)\,dt,
\]
which identifies $P(\rho)/2$ with the flux of the relative stress
across the slice $\{\rho=\mathrm{const}\}$ in the
$\nabla t$-direction.

Next, we multiply the identity
\eqref{eq:stress-rho} by $t$ and integrate in $t$ from $0$ to $\rho$.
Using $t\omega_\rho=-\rho\omega_t$ and integrating by parts in $t$, we obtain
\[
\frac{d}{d\rho}
\int_0^\rho t\omega(2S_\rho T_\rho+e_{\mathrm{rel}})\,dt
=\frac{\rho P'-P}{2}.
\]
Integrating in $\rho$ from $0$ to $\rho$ then yields
\[
\int_0^\rho t\omega(2S_\rho T_\rho+e_{\mathrm{rel}})\,dt
=\frac{\rho P(\rho)}2-\int_0^\rho P(\rho')\,d\rho'.
\]

Now set
\[
F(\rho)=\rho^{-3}\int_0^\rho P(\rho')\,d\rho'.
\]
Combining the preceding identity with
\eqref{eq:Pprime} and using $P=(\rho^3F)'$, we obtain
\begin{equation}\label{eq:stress-convexity}
	\begin{aligned}
	&F''-\frac{n-8}{\rho}F'-\frac{n-8}{\rho^2}F\\
	&\qquad=
	\frac{2(n-2)}{\rho^5}\int_0^\rho t\omega
	\left[
	-2S_\rho T_\rho
	+\frac{t^2}{\rho^2-t^2}e_{\mathrm{rel}}
	\right]dt>0,
	\end{aligned}
\end{equation}
where the positivity of the right-hand side follows from
$S_\rho T_\rho<0$ and $e_{\mathrm{rel}}>0$.

Since $P(0)=0$ and $P'>0$, we have $F>0$ for $\rho>0$.
Smoothness and $O(m)\times O(m)$-invariance, together with
$S(0)=0$, yield $S=O(|x|^2)$ and
$\nabla S,\nabla T=O(|x|)$ near the origin. Consequently, as $\rho\downarrow0$,
\[
P(\rho)=O(\rho^{n+1}),\qquad
F(\rho)=O(\rho^{n-1}),\qquad
F'(\rho)=O(\rho^{n-2}).
\]
Since $n\ge8$, \eqref{eq:stress-convexity} implies that
\[
\bigl(\rho^{8-n}F'\bigr)'
>(n-8)\rho^{6-n}F\ge0.
\]
Combined with
$\lim_{\rho\downarrow0}\rho^{8-n}F'(\rho)=0$,
this gives $F'>0$. Returning to
\eqref{eq:stress-convexity}, we further obtain $F''>0$,
and hence $F(\rho)\ge c\rho$ for some constant $c>0$ and all sufficiently large $\rho$.

On the other hand, the boundedness of $S$ and $\nabla S$, together
with Proposition~\ref{prop:pointwise-estimates}\textup{(i)} and \textup{(ii)}, yields
\[
0<e_{\mathrm{rel}}(\rho,t)
\le
	\begin{cases}
	C\rho^{-\beta}e^{-\mu t},
	&0<t<\kappa\rho,\\
	Ce^{-c\rho},
	&\kappa\rho\le t<\rho
	\end{cases}
\]
for all sufficiently large $\rho$.
Using \eqref{eq:Pprime} and
$(\rho^2-t^2)^{m-2}\le\rho^{n-4}$, we obtain
\[
	\begin{aligned}
	0<P'(\rho)
	&\le C\rho^{n-4-\beta}
	\int_0^{\kappa\rho}t e^{-\mu t}\,dt
	+C\rho^{n-2}e^{-c\rho}\\
	&\le C\rho^{\alpha-1},
	\end{aligned}
\]
where $\alpha+\beta=n-3$.
Integrating this estimate yields $P(\rho)\le C\rho^\alpha$,
and consequently
\[
F(\rho)=\frac{1}{\rho^3}
\int_0^\rho P(\rho')\,d\rho'
\le C\rho^{\alpha-2}
\]
for all sufficiently large $\rho$.
Since $\alpha\le2$, this contradicts the lower bound
$F(\rho)\ge c\rho$. We conclude that $T\equiv0$.
\end{proof}

\begin{proof}[Proof of Theorem~\ref{thm:main}]
Lemma~\ref{lem:rigidity} yields
$
U_0^+=U_0^-=:U_0.
$
By Lemma~\ref{lem:envelope}, every $V\in\M_0$ lies between these two extremal elements, and hence $\M_0=\{U_0\}$.
Moreover, the reflection identity implies $U_0(r,s)=-U_0(s,r)$, while
Lemma~\ref{lem:coordinate-signs} gives
$U_{0,r}>0$ and $U_{0,s}<0$. Since $U_0(s,s)=0$, it follows that
$U_0>0$ for $r>s$ and $U_0<0$ for $r<s$.
Thus $U_0$ is a saddle-shaped global minimizer.

Finally, by the uniqueness
of saddle-shaped solutions established in \cite[Theorem~1.2]{Cabre2012},
$U_0$ coincides with the saddle-shaped solution constructed by
Cabr\'e and Terra in \cite{CabreTerra2009}.
\end{proof}

\appendix

\section{Auxiliary estimates}
\label{app:identities}

To obtain the estimates in Section~\ref{sec:refined}, we multiply $u-H$ by a cutoff supported near the cone, extend by zero in the normal variable, and decompose the result, at each fixed $\rho$, into $\xi(\rho)H'(t)$ and a remainder orthogonal to $H'$ in $L^2(\mathbb{R},dt)$. Since $H(t+a)=H(t)+aH'(t)+O(a^2)$, the first component represents a normal displacement of the interface to first order. Subtracting the  correction $(n-2)\rho^{-2}\eta(t)$ from the orthogonal component cancels the leading geometric error and leaves a smaller remainder, which we estimate using the spectral gap of the one-dimensional linearized operator. Projection onto $H'$ yields an inhomogeneous radial Jacobi equation for $\xi$ with controlled decaying forcing. 

We begin with the formulas for the gradient and Laplacian in $(\rho,t)$ coordinates, the spectral properties of the one-dimensional linearized operator, and the construction of $\eta$. These estimates determine the coefficient $A$ of $\rho^{-\alpha}H'$ in the expansion of $u-H$. For two solutions with the same coefficient $A$, this term cancels in their difference. The equation for $u-v$, together with the scalar ODE estimates in Section~\ref{sec:refined}, then improves its decay from $\rho^{-\alpha}$ to $\rho^{-\beta}$, where $\beta>\alpha$. Finally, we differentiate the Allen--Cahn equation and the system satisfied by the odd and even parts of a solution under $r\leftrightarrow s$. The resulting equations are used in the main text to prove the required inequalities for the derivatives of the extremal minimizers and their odd and even parts.

Recall from Section~\ref{sec:intro} that $n=2m$ with $m\ge4$, and
\[
r=|x'|,\quad s=|x''|,\quad
\rho=\frac{r+s}{\sqrt2},\quad t=\frac{r-s}{\sqrt2}.
\]
The numbers $\alpha<\beta$ are the roots of $k^2-(n-3)k+n-2=0$. Thus, $1<\alpha\le2$, $\beta>\alpha$, $\alpha+\beta=n-3$ and $\alpha\beta=n-2$. We fix
\begin{equation}\label{apeq:epsilon}
	0<\epsilon_*<\min\left\{
	\frac{\alpha-1}{4},\frac{\beta-\alpha}{4},\frac12\right\}.
\end{equation}
Unless stated otherwise, functions are $O(m)\times O(m)$-invariant. 

\subsection{Coordinates and basic operators}
\label{app:coordinates-subsection}

Write
\[
\omega=(\rho^2-t^2)^{m-1},\quad
A=\frac{(n-2)\rho}{\rho^2-t^2},\quad
B=-\frac{(n-2)t}{\rho^2-t^2}.
\]
A simple direct computation using the chain rule gives the following identities; see also \cite[(1.12), (1.14)--(1.15)]{CabreTerra2010}, whose variables $(s,t,y,z)$ correspond to our $(r,s,\rho,t)$.

\begin{lemma}\label{lem:coordinates}
	On $\{r,s>0\}$,
	\begin{equation}\label{apeq:gradient-identities}
		|\nabla\rho|=|\nabla t|=1,\qquad
		\nabla\rho\cdot\nabla t=0,\qquad
		\Delta\rho=A,\qquad\Delta t=B.
	\end{equation}
 Since $\partial_r=\frac{\partial_\rho+\partial_t}{\sqrt2}$ and $\partial_s=\frac{\partial_\rho-\partial_t}{\sqrt2}$, we have
	\begin{align}
		\Delta&=\partial_{rr}+\partial_{ss}
		+\frac{m-1}{r}\partial_r+\frac{m-1}{s}\partial_s
		\notag\\
		&=\partial_{\rho\rho}+\partial_{tt}
		+A\partial_\rho+B\partial_t,
		\label{apeq:laplacian}
	\end{align}
	and for an invariant function $f$, we have 
	\begin{equation}
		\Delta f=\omega^{-1}
		\bigl[\partial_\rho(\omega f_\rho)+\partial_t(\omega f_t)\bigr].
	\end{equation}
\end{lemma}

Now consider the following operators.
\begin{equation}\label{apeq:operators}
	\begin{gathered}
		L_H=-\partial_t^2+3H^2-1,\qquad
		\mathcal L_U=-\Delta+3U^2-1,\qquad
		\mathcal L_H=-\Delta+3H(t)^2-1,\\
		\mathcal L_0=-\partial_{\rho\rho}
		-\frac{n-2}{\rho}\partial_\rho+L_H,
		\qquad
		J\xi=\xi''+\frac{n-2}{\rho}\xi'+\frac{n-2}{\rho^2}\xi.
	\end{gathered}
\end{equation}
In particular,
\begin{equation}\label{apeq:model-difference}
	\mathcal L_H
	=\mathcal L_0-\left(A-\frac{n-2}{\rho}\right)\partial_\rho
	-B\partial_t,
	\qquad
	A-\frac{n-2}{\rho}=\frac{(n-2)t^2}{\rho(\rho^2-t^2)}.
\end{equation}

\subsection{Linearizations}
\label{acapp:heteroclinic-subsection}


\begin{lemma}\label{lem:normal}
	The operator $L_H:H^2(\mathbb R)\to L^2(\mathbb R)$ has
	$\ker L_H=\operatorname{span}\{h\}$, and
	\begin{equation}\label{apeq:spectral-gap}
		\int_{\mathbb R}\bigl(|f'|^2+(3H^2-1)f^2\bigr)\,dt
		\ge\frac32\int_{\mathbb R}f^2\,dt
	\end{equation}
	if $f\in H^1(\mathbb R)$ is such that $\int_{\mathbb R}fh\,dt=0$.
	
	There is a unique odd $\eta\in H^2(\mathbb R)$ satisfying
	\begin{equation}\label{apeq:eta-equation}
		L_H\eta=-th,
	\end{equation}
	which is smooth and, for every integer $j\ge0$ and $0<\mu<\sqrt2$, satisfies
	\begin{equation}\label{apeq:eta-decay-signs}
		|\eta^{(j)}(t)|\le C_{j,\mu}e^{-\mu|t|},\qquad
		\eta(t)<0\ (t>0),\qquad \eta'(0)<0,
	\end{equation}
	Moreover,
	\begin{equation}\label{apeq:normal-identities}
		\int_\R th'h\,dt=-\frac\sigma2,\qquad
		\int_\R H\eta h^2\,dt=-\frac\sigma{12},\qquad
		\int_\R Hh^3\,dt=0.
	\end{equation}
\end{lemma}

\begin{proof}
	For the kernel and spectral gap of $L_H$, see \cite[Lemma~8.1 and (8.3)]{PW}. The correction $\eta$ is the function $\psi_0$ defined immediately before (3.8) in \cite[Section~3.2]{DPKW2013} (with $w=H$). Its oddness and decay, as well as the signs $\eta(t)<0$ for $t>0$ and $\eta'(0)<0$, follow from the integral formula given there.
	The derivative bounds follow by differentiation. The first and third identities in \eqref{apeq:normal-identities} follow by integration by parts and parity; the second is the identity in \cite[proof of Lemma~9.1, p.~111]{DPKW2013}, with $f(s)=s-s^3$ and $w=H$.
\end{proof}

\subsection{Cut-off, projection, and orthogonal estimates}
\label{acapp:linear-subsection}
We work in the $(\rho,t)$-plane, suppressing the angular variables by $O(m)\times O(m)$-invariance. The region $\{r,s>0\}$ corresponds to $\rho>0$ and $|t|<\rho$. For fixed $R>0$ and $0<\kappa<1$, we write $\mathcal T_{\kappa,R}=\{(\rho,t):\rho>R,\ |t|<\kappa\rho\}$. After multiplication by a cutoff supported in this region, we can extend functions on $\mathcal T_{\kappa,R}$ by zero (in $t$) and regard them as functions on the larger domain $D_R=(R,\infty)\times\mathbb R$.

For a domain $D$ in the $(\rho,t)$-plane, $j\in\mathbb N\cup\{0\}$, $0<\gamma<1$, $\ell\in\mathbb R$ and $\mu>0$, set
\begin{equation}\label{apeq:weighted-norm}
	\|F\|_{C^{j,\gamma}_{\ell,\mu}(D)}
	=\sup_{X_0=(\rho_0,t_0)\in D}
	\rho_0^\ell e^{\mu|t_0|}
	\|F\|_{C^{j,\gamma}(B_1(X_0)\cap D)}.
\end{equation}
Moreover, we set 
\[
\|F\|_{C^0_{\ell,\mu}(D)}=\sup_{(\rho,t)\in D}\rho^\ell e^{\mu|t|}|F(\rho,t)|. 
\]
We write $F=O_\mu(\rho^{-\ell})$ in $D$ if $|F(\rho,t)|\le C\rho^{-\ell}e^{-\mu|t|}$ throughout $D$. We write $F=O_{\mu,\gamma}(\rho^{-\ell})$ if, more strongly, $\|F\|_{C^{0,\gamma}(B_1(X_0)\cap D)}\le C\rho_0^{-\ell}e^{-\mu|t_0|}$ for every $X_0=(\rho_0,t_0)\in D$. 

For every fixed $\rho$ let
\begin{equation}\label{acapp:projection-definition}
	\Pi F=F-\frac{h}{\sigma}\int_{\mathbb R}Fh\,dt
\end{equation}
be the orthogonal projection onto $h^\perp$ in $L^2(\mathbb R,dt)$.

\begin{lemma}\label{lem:localization}
	Fix $0<\kappa<\kappa_1<1/4$, and choose an even cutoff function $\chi\in C_c^\infty((-\kappa_1,\kappa_1))$ such that $0\le\chi\le1$ and $\chi=1$ on $[-\kappa,\kappa]$.

	Set $\vartheta(\rho,t)=\chi(t/\rho)$. For a smooth function $v$ on $\mathcal T_{\kappa_1,R}$, extend $\vartheta v$ by zero to $D_R$.
	Then
	\begin{equation}\label{apeq:cutoff-identity1}
			\mathcal L_0(\vartheta v) = \vartheta\mathcal L_0v+\mathcal C_\chi[v],
	\end{equation}
	where 
	\begin{equation}\label{apeq:cutoff-identity2}
		\mathcal C_\chi[v] =-2\vartheta_\rho v_\rho-2\vartheta_t v_t
			-\left(\vartheta_{\rho\rho}
			+\frac{n-2}{\rho}\vartheta_\rho+\vartheta_{tt}\right)v.
	\end{equation}
	If $|v|+|v_\rho|+|v_t|\le C\rho^N e^{-\mu_0|t|}$, then for every
	$0<\mu<\min\{\mu_0,\sqrt2\}$ there is $c>0$ such that
	\begin{equation}\label{apeq:cutoff-smallness}
		|\mathcal C_\chi[v]|+|\Pi\mathcal C_\chi[v]|
		\le Ce^{-c\rho}e^{-\mu|t|}\quad\text{in }D_R.
	\end{equation}
	If $v,v_\rho,v_t\in C^{0,\gamma}_{-N,\mu_0}$, the same estimate holds in the weighted local $C^{0,\gamma}$ norms, after decreasing $c$ if necessary. 
\end{lemma}

\begin{proof}
    The product rule gives \eqref{apeq:cutoff-identity1} and \eqref{apeq:cutoff-identity2}. Any positive-order derivative of $\vartheta$ is supported where $|t|\ge\kappa\rho$ and is bounded by a polynomial in $\rho$. A portion of the factor $e^{-(\mu_0-\mu)|t|}$ absorbs this polynomial growth, while the remainder produces an exponentially small bound of $O(e^{-c\rho})$. This proves the estimate for $\mathcal C_\chi[v]$ in \eqref{apeq:cutoff-smallness}. Because its projection coefficient is also $O(e^{-c\rho})$ and $h$ decays at a rate of $\sqrt2$, the bound for $\Pi\mathcal C_\chi[v]$ directly follows. The product estimate then gives the H\"older version. The final assertion holds because $\vartheta$ and $h$ are even by construction.
\end{proof}

When extending expressions from $\mathcal T_{\kappa,R}$ to $D_R$, we first multiply them by the prescribed cutoff function and then extend the result by zero. This procedure ensures that the coefficients remain smooth across $t=\pm\rho$.

\begin{lemma}\label{lem:projection}
Suppose that $\mathcal L_0w=F$ in $D_R$, and that $w$, its coordinate derivatives up to order two, and $F$ all decay exponentially in $|t|$, locally uniformly in $\rho$. More precisely, for every compact interval $I\subset(R,\infty)$, there exist constants $C_I,\mu_I>0$ such that
\[
	\sum_{\substack{j,k\ge0\\j+k\le2}}|\partial_\rho^j\partial_t^k w(\rho,t)|+|F(\rho,t)|
	\le C_I e^{-\mu_I|t|}
	\qquad\text{for all }(\rho,t)\in I\times\mathbb R.
\]
Let us define
\[
	\zeta(\rho)=\frac1\sigma\int_{\mathbb R}wh\,dt,
	\quad p=w-\zeta h,
	\quad \text{and} \quad f(\rho)=\frac1\sigma\int_{\mathbb R}Fh\,dt.
\]
Then, we have the identities
\begin{equation}\label{apeq:projection-identities}
	J\zeta=\frac{n-2}{\rho^2}\zeta-f,
	\qquad \text{and}\qquad \mathcal L_0p=\Pi F.
\end{equation}
Furthermore, for any function $\tau\in C^2((R,\infty))$, it holds that
\begin{equation}\label{apeq:corrected-projection}
	\mathcal L_0(p-\tau\eta)
	=\Pi F+\left(\tau''+\frac{n-2}{\rho}\tau'\right)\eta+\tau th.
\end{equation}
\end{lemma}

\begin{proof}
	Since $L_Hh=0$, projecting onto $h$ and integrating by parts in $t$ gives
	$-\zeta''-(n-2)\rho^{-1}\zeta'=f$.
	This is the first identity in \eqref{apeq:projection-identities}. Orthogonal projection gives the second. Finally, \eqref{apeq:corrected-projection} follows from $L_H\eta=-th$ and the product rule.
\end{proof}

\begin{lemma}[Estimate on the orthogonal component]\label{lem:orthogonal}
	Let $R\ge1$ and suppose $\mathcal L_0p=G$ in $D_R$, with $\int_{\mathbb R}p(\rho,t)h(t)\,dt=0$ for every $\rho>R$. Assume $p$ is smooth and
	\[
	\sup_{\rho\in I}\sum_{j+\ell\le2}
	\|\partial_\rho^j\partial_t^\ell p(\rho,\cdot)\|_{L^2(\mathbb R)}
	<\infty
	\]
	for every compact $I\subset(R,\infty)$. Then $y(\rho)=\|p(\rho,\cdot)\|_{L^2(\mathbb R)}$ satisfies
	\begin{equation}\label{apeq:L2-inequality}
		-y''-\frac{n-2}{\rho}y'+\frac32y
		\le\|G(\rho,\cdot)\|_{L^2(\mathbb R)}
	\end{equation}
	in the distributional sense.
	
	Suppose in addition that $s,k,\mu_0,\mu>0$, $0<\gamma<1$, and
	\[
	M_0=\|p\|_{C^0_{s,\mu_0}(D_R)}<\infty,
	\qquad M_1=\|G\|_{C^{0,\gamma}_{k,\mu}(D_R)}<\infty.
	\]
	For every $0<\mu'<\min\{\mu,\mu_0,\sqrt2\}$ there are
	$R'>R+4$ and $C$, depending only on
	$n,R,s,k,\mu_0,\mu,\mu',\gamma$, such that
	\begin{equation}\label{apeq:orthogonal-estimate}
		\|p\|_{C^{2,\gamma}_{k,\mu'}(D_{R'})}
		\le C(M_0+M_1).
	\end{equation}
\end{lemma}

\begin{proof}
	We first use the spectral gap to control the $L^2$ norm. Where $y>0$, we have $yy'=\langle p,p_\rho\rangle$ and $yy''=\|p_\rho\|_2^2+\langle p,p_{\rho\rho}\rangle-(y')^2\ge\langle p,p_{\rho\rho}\rangle$.
	Then \eqref{apeq:spectral-gap} implies
	\[
	y\left(-y''-\frac{n-2}{\rho}y'\right)
	\le-\langle p,L_Hp\rangle+\langle p,G\rangle
	\le-\frac32y^2+y\|G(\rho,\cdot)\|_2.
	\]
	Applying the same calculation to $(y^2+\varepsilon^2)^{1/2}$ and letting $\varepsilon\downarrow0$ proves \eqref{apeq:L2-inequality} across the zero set of $y$.
	
	For the quantitative estimate, we note that $y(\rho)\le C M_0\rho^{-s}\to0$ and $\|G(\rho,\cdot)\|_2\le C M_1\rho^{-k}$. Provided $R_1>R+1$ is sufficiently large (depending only on the fixed parameters), it follows that
\[
	\left(-\frac{d^2}{d\rho^2}-\frac{n-2}{\rho}\frac d{d\rho}
	+\frac32\right)\rho^{-k}
	=\left(\frac32+k(n-3-k)\rho^{-2}\right)\rho^{-k}
	\ge\rho^{-k}.
\]
    We then select a multiple $C(M_0+M_1)\rho^{-k}$ that simultaneously dominates $y$ at $\rho=R_1$. By the comparison principle applied to \eqref{apeq:L2-inequality}, together with the fact that the solution vanishes at infinity, we obtain
\begin{equation}\label{apeq:L2-decay}
	y(\rho)\le C(M_0+M_1)\rho^{-k}
\end{equation}  
    on the domain $\rho\ge R_1$.
    
	Set $M=M_0+M_1$, and let $B_a(X_0)$ denote the Euclidean ball of radius $a$ centred at $X_0=(\rho_0,t_0)$ in the $(\rho,t)$-plane. For $\rho_0\ge R_1+2$, the preceding bound on $y$ and the definition of $M_1$ give
	\[
	\|p\|_{L^2(B_1(X_0))}+\|G\|_{L^\infty(B_1(X_0))}\le CM\rho_0^{-k}.
	\]
	The coefficients of $\mathcal L_0$ are uniformly bounded on these balls. Hence the interior local boundedness estimate, applied to $p$ and $-p$, yields
	\begin{equation}\label{apeq:pointwise-polynomial}
		|p(\rho,t)|\le CM\rho^{-k}
	\end{equation}
	for $\rho\ge R_1+2,\ t\in\mathbb R$.
	
	To recover exponential decay in $t$, write $V(t)=3H(t)^2-1$ and choose $t_*>0$ and $R_2\ge R_1+2$ sufficiently large that
	\[
	V(t)\ge c_0,\qquad V(t)-(\mu')^2+\frac{k(n-3-k)}{\rho^2}\ge c_0
	\]
	for $\rho\ge R_2$ and $|t|\ge t_*$, with $c_0>0$. This is possible because $V(t)\to2$ and $\mu'<\sqrt2$. On each domain $\Omega_\pm=\{(\rho,t):\rho>R_2,\ \pm t>t_*\}$, the function $b(\rho,t)=\rho^{-k}e^{-\mu'|t|}$ satisfies
	\[
	\mathcal L_0b=\left(V(t)-(\mu')^2+\frac{k(n-3-k)}{\rho^2}\right)b\ge c_0b.
	\]
	Since $\mu,\mu_0>\mu'$, a sufficiently large multiple $\Phi=CMb$ satisfies $\mathcal L_0\Phi\ge|G|$ in $\Omega_\pm$ and $\Phi\ge|p|$ on $\partial\Omega_\pm$. For the boundary comparison, use \eqref{apeq:pointwise-polynomial} on $t=\pm t_*$ and the original bound $|p(R_2,t)|\le M_0R_2^{-s}e^{-\mu_0|t|}$ on $\rho=R_2$. Both $p$ and $\Phi$ tend to zero as $\rho+|t|\to\infty$. The maximum principle applied to $p-\Phi$ and $-p-\Phi$ therefore gives $|p|\le\Phi$ in $\Omega_\pm$. On $|t|\le t_*$, the same estimate follows from \eqref{apeq:pointwise-polynomial}, after enlarging the constant. Thus, we obtain exponential decay in $t$
	\[
	|p(\rho,t)|\le CM\rho^{-k}e^{-\mu'|t|}
	\]
	for $\rho\ge R_2,\ t\in\mathbb R$.
	
	Finally, set $R'=R_2+3$. Since the coefficients of $\mathcal L_0$ have uniformly bounded $C^{0,\gamma}$ norms, and the weights are uniformly comparable on balls of radius two, the interior Schauder estimate  yields
    \[
        \|p\|_{C^{2,\gamma}(B_1(X_0))}\le CM\rho_0^{-k}e^{-\mu'|t_0|}
    \]
	for $X_0=(\rho_0,t_0)\in D_{R^{\prime}}$. Taking the weighted supremum over $X_0\in D_{R^{\prime}}$ proves \eqref{apeq:orthogonal-estimate}, with all constants and radii depending only on the fixed parameters.

\end{proof}

When applying these estimates repeatedly, we fix nested tubular domains in advance and increase the exterior radius whenever necessary. Since $r,s\asymp\rho$ on any such fixed tube, the derived coordinate bounds immediately yield the corresponding Euclidean estimates.

\subsection{Projection of the coarse profile}
\label{acapp:coarse-subsection}

Let $u$ solve $-\Delta u=u-u^3$ and put $q=u-H$. Since
$\Delta H=H^3-H+Bh$, we have
\begin{equation}\label{apeq:profile-equation}
	\mathcal L_0q
	=Bh+\left(A-\frac{n-2}{\rho}\right)q_\rho+Bq_t-3Hq^2-q^3.
\end{equation}
The geometric term $Bh$ on the right-hand side satisfies $Bh=-(n-2)\rho^{-2}th+O_\mu(\rho^{-4})$ on every fixed tube $\mathcal T_{\kappa,R}$ with $\kappa<1$, for $0<\mu<\sqrt2$. The correction $(n-2)\rho^{-2}\eta$, where $L_H\eta=-th$, is chosen to cancel the displayed term of order $\rho^{-2}$; compare \cite[Section~2.2, proof of Proposition~11]{LWW}.

\begin{lemma}\label{lem:coarse}
Fix $0<\kappa<\kappa_1<1/4$. Let $u$ be an invariant entire solution of $-\Delta u=u-u^3$ satisfying
\begin{equation}\label{apeq:coarse-hypothesis}
	|D^j(u-H)|\le M\rho^{-\alpha}e^{-\mu_0|t|}
	\quad(0\le j\le2)\quad\text{in }\mathcal T_{\kappa_1,R}
\end{equation}
for given constants $M,\mu_0>0$ and $R\ge1$. Let $\vartheta$ be the cutoff function introduced in Lemma~\ref{lem:localization}, and define on $D_R$ the quantities
\begin{equation}\label{apeq:coarse-decomposition}
	\widetilde q=\vartheta q,\quad
	\xi=\frac1\sigma\int_{\mathbb R}\widetilde qh\,dt,\quad
	p=\widetilde q-\xi h-(n-2)\rho^{-2}\eta.
\end{equation}
It then follows that $p\perp h$. By increasing $R$ and decreasing the exponential rate if needed, we secure the bound $|\xi|\le C\rho^{-\alpha}$ as well as
\begin{equation}\label{apeq:coarse-orthogonal}
	\|p(\rho,\cdot)\|_2+\sum_{j+\ell\le2}\sup_{t\in\mathbb R}
	e^{\mu|t|}|\partial_\rho^j\partial_t^\ell p(\rho,t)|
	\le C\rho^{-2\alpha}.
\end{equation}
Furthermore,
\begin{equation}\label{apeq:coarse-Jacobi}
	|J\xi|\le C\bigl(\rho^{-3\alpha}+\rho^{-\alpha-3}\bigr)
	\le C\rho^{-\alpha-2-2\epsilon_*}.
\end{equation}
On the restricted tube $\mathcal T_{\kappa,R}$, the derivative estimate also holds for Euclidean derivatives. The implicit constants, along with the adjusted radius and exponential rate, depend only on the parameters $n$, $R$, $\kappa$, $\kappa_1$, $\chi$, $M$, $\mu_0$, and the fixed exponents.
\end{lemma}

\begin{proof}
	Fix $\gamma=1/2$. The bounds in \eqref{apeq:coarse-hypothesis} give local $C^{0,\gamma}$ bounds for $q,q_\rho,q_t$ with the same weights. Write
	\[
	\mathcal L_0\widetilde q=F
	\coloneqq\vartheta\left[
	Bh+\left(A-\frac{n-2}{\rho}\right)q_\rho
	+Bq_t-3Hq^2-q^3\right]+\mathcal C_\chi[q].
	\]
	The cutoff error $\mathcal C_\chi[q]$ is exponentially small by Lemma~\ref{lem:localization}. The initial estimates give $|\xi|+\|p(\rho,\cdot)\|_2=O(\rho^{-\alpha})$ and a similar pointwise bound for $p$, since $\alpha\le2$.
	
	By \eqref{apeq:corrected-projection},
	\begin{equation}\label{apeq:coarse-perp-equation}
		\mathcal L_0p
		=\Pi F+(n-2)\rho^{-2}th
		-2(n-2)(n-5)\rho^{-4}\eta\eqqcolon F_\perp.
	\end{equation}
	The term $\vartheta Bh$ is odd, so its projection onto $h$ vanishes. Its leading part cancels the added $th$ term:
	\begin{equation}\label{apeq:curvature-cancellation}
		-\frac{(n-2)\vartheta t}{\rho^2-t^2}h+(n-2)\rho^{-2}th
		=-\frac{(n-2)\vartheta t^3}{\rho^2(\rho^2-t^2)}h
		+(n-2)\rho^{-2}(1-\vartheta)th.
	\end{equation}
	The first term is $O_\mu(\rho^{-4})$ and the second is exponentially small in $\rho$. Estimating the remaining terms in \eqref{apeq:coarse-perp-equation}, including their projections, gives
	\[
	F_\perp
	=O_{\mu,\gamma}
	(\rho^{-2\alpha}+\rho^{-\alpha-2}+\rho^{-4})
	=O_{\mu,\gamma}(\rho^{-2\alpha}).
	\]
	Here we used $\alpha\le2$ and absorbed powers of $|t|$ by decreasing $\mu$. Lemma~\ref{lem:orthogonal} then gives \eqref{apeq:coarse-orthogonal}.
	
	We now compute $\sigma^{-1}\int Fh$. There are two leading contributions,
	one from the normal drift and one from the quadratic interaction with
	$\eta$. Each equals $(n-2)\xi/(2\rho^2)$.
	
	The term $\vartheta Bh$, which is independent of $q$, has zero projection onto $h$ by parity. The term $\vartheta(A-(n-2)/\rho)q_\rho$ contributes $O(\rho^{-\alpha-3})$ to the scalar projection, since $A-(n-2)/\rho=(n-2)t^2/[\rho(\rho^2-t^2)]$ and $q_\rho=O_\mu(\rho^{-\alpha})$ in the tube. For the term $\vartheta Bq_t$, involving the derivative of $q$ in the direction normal to the cone, we use $\vartheta q_t=\partial_t\widetilde q-\vartheta_tq$ to compute its projection onto $h$.
	\begin{align}
		\int_{\mathbb R}\vartheta Bq_t h\,dt
		&=-\frac{n-2}{\rho^2}\int_{\mathbb R}t\widetilde q_t h\,dt
		+O(\rho^{-\alpha-4})\notag\\
		&=\frac{(n-2)\sigma}{2\rho^2}\xi
		+O(\rho^{-2\alpha-2}+\rho^{-\alpha-4}).
		\label{apeq:normal-drift-projection}
	\end{align}
	The second equality uses \eqref{apeq:coarse-decomposition}, $\int th'h=-\sigma/2$, $\int t\eta'h=0$, and \eqref{apeq:coarse-orthogonal}.
	
	Replacing $\vartheta q^j$ by $\widetilde q^{\,j}$ for $j=2,3$
	changes the projected integrals only by an exponentially small term.
	Expanding \eqref{apeq:coarse-decomposition} gives
	\begin{align}
		-\int_{\mathbb R}\vartheta(3Hq^2+q^3)h\,dt
		&=-6(n-2)\rho^{-2}\xi
		\int_{\mathbb R}H\eta h^2\,dt
		+O(\rho^{-3\alpha})\notag\\
		&=\frac{(n-2)\sigma}{2\rho^2}\xi+O(\rho^{-3\alpha}).
		\label{apeq:nonlinear-projection}
	\end{align}
	The terms $-3\xi^2Hh^2$ and $-3(n-2)^2\rho^{-4}H\eta^2$ in the expansion of $-3H\widetilde q^{\,2}$ have zero projection onto $h$, since $Hh^3$ and $H\eta^2h$ are odd functions of $t$. Terms involving $p$, and the cubic terms, have the stated remainder bound because $p=O_\mu(\rho^{-2\alpha})$ and $\alpha\le2$.
	
	Adding \eqref{apeq:normal-drift-projection} and \eqref{apeq:nonlinear-projection}, we find
	\begin{equation}\label{acapp:coarse-projected-forcing}
		\frac1\sigma\int_{\mathbb R}Fh\,dt
		=\frac{n-2}{\rho^2}\xi
		+O(\rho^{-3\alpha}+\rho^{-\alpha-3}).
	\end{equation}
	The first identity in \eqref{apeq:projection-identities} now gives \eqref{apeq:coarse-Jacobi}. Its final bound follows from \eqref{apeq:epsilon}.
\end{proof}

The hypotheses of this lemma are provided by Proposition~\ref{prop:trapping-bounds} for solutions lying between fixed bounding minimizers, and the associated constants are uniform across all such solutions. Applying the scalar ODE estimate of Lemma~\ref{lem:euler} to \eqref{apeq:coarse-Jacobi} yields the expansion given in Lemma~\ref{lem:coarse-expansion}. While the correction term $\eta$ is orthogonal to $h$, its interaction with $\xi h$ in \eqref{apeq:nonlinear-projection} is crucial: combined with the normal drift, it generates the complete Jacobi potential.

\subsection{Projection of the difference equation}
\label{acapp:difference-subsection}

For the difference $w=u-v$, we first introduce the variable $z=\Lambda w$, where
\[
	\Lambda=\left(1-\frac{t^2}{\rho^2}\right)^{(n-2)/4}.
\]
This transformation eliminates the first-order normal derivative term $-Bw_t$ and replaces the radial coefficient $A$ with $(n-2)/\rho$. The resulting equation involves the operator $\mathcal L_0$ and an additional zeroth-order coefficient, which simplifies its projection onto $h$ and its orthogonal complement.

\begin{lemma}\label{lem:conjugation}
	For an arbitrary potential $Q$, we have on the domain $|t|<\rho$ that
	\begin{equation}\label{apeq:conjugation-identity}
		\Lambda(-\Delta+Q)\Lambda^{-1}
		=-\partial_{\rho\rho}-\frac{n-2}{\rho}\partial_\rho
		-\partial_{tt}+Q+V_0,
	\end{equation}
	where
	\begin{equation}\label{apeq:geometric-potential}
		V_0=\frac{(n-2)(n-6)(\rho^2+t^2)}{4(\rho^2-t^2)^2}
		-\frac{(n-2)(n-4)}{4\rho^2}.
	\end{equation}
	Moreover, within a fixed tubular region, we have the asymptotic expansion
	\begin{equation}\label{apeq:geometric-expansion}
		V_0=-\frac{n-2}{2\rho^2}+O(t^2\rho^{-4}).
	\end{equation}
\end{lemma}

\begin{proof}
Let us set $g=\log\Lambda$. Differentiating the relation $2g=\log\omega-(n-2)\log\rho$ yields $2g_t=B$ and $2g_\rho=A-(n-2)/\rho$. By the product rule, the coefficients of $\partial_t$ and $\partial_\rho$ in the conjugated operator $\Lambda(-\Delta)\Lambda^{-1}$ are therefore $2g_t-B=0$ and $2g_\rho-A=-(n-2)/\rho$, respectively. The zeroth-order term is given by
\begin{align*}
	V_0
	&=g_{\rho\rho}+g_{tt}+Ag_\rho+Bg_t-g_\rho^2-g_t^2\\
	&=\frac12(A_\rho+B_t)+\frac14(A^2+B^2)-\frac{(n-2)(n-4)}{4\rho^2}.
\end{align*}
Substituting the explicit expressions for $A$ and $B$, we obtain \eqref{apeq:geometric-potential}. Since multiplication by $Q$ is invariant under conjugation, this establishes the stated operator identity. Finally, expanding the term $(\rho^2+t^2)/(\rho^2-t^2)^2=\rho^{-2}+O(t^2\rho^{-4})$ in the region $|t|\le\kappa_1\rho$ yields \eqref{apeq:geometric-expansion}.
\end{proof}

\begin{lemma}\label{lem:difference}
	Fix $0<\kappa<\kappa_1<\kappa_2<\kappa_3<1/4$. Let $u$ and $v$ be invariant entire solutions satisfying \eqref{apeq:coarse-hypothesis} on $\mathcal T_{\kappa_3,R}$; by enlarging the constants if necessary, we may assume they share the same constants. Choose the cutoff $\vartheta$ as in Lemma~\ref{lem:localization}, supported in the $\kappa_1$-tube and identically equal to one on the $\kappa$-tube. Setting $w=u-v$, we define
	\begin{equation}\label{apeq:difference-decomposition}
		\widetilde w=\vartheta\Lambda w,\quad
		\xi=\frac1\sigma\int_{\mathbb R}\widetilde wh\,dt,
		\quad \text{and} \quad p=\widetilde w-\xi h.
	\end{equation}
	Next, let $\widehat\chi$ be a smooth even function such that $0\le\widehat\chi\le1$, supported in $(-\kappa_2,\kappa_2)$ and identically one on a neighbourhood of $\operatorname{supp}\chi$. We set $\widehat\vartheta(\rho,t)=\widehat\chi(t/\rho)$, define the potential
	\[
	Q_\Lambda=\widehat\vartheta
	(u^2+uv+v^2-3H^2+V_0),
	\]
	and extend it by zero to the domain $D_R$. Then, provided $R$ is chosen sufficiently large, we have
	\begin{equation}\label{apeq:difference-equation}
		\mathcal L_0\widetilde w
		=-Q_\Lambda\widetilde w+\mathcal C_\chi[\Lambda w],
		\quad \|Q_\Lambda(\rho,\cdot)\|_\infty\le C\rho^{-\alpha},
	\end{equation}
	and
	\begin{equation}\label{apeq:averaged-potential}
		\frac1\sigma\int_{\mathbb R}Q_\Lambda h^2\,dt
		=-\frac{n-2}{\rho^2}
		+O(\rho^{-2\alpha}+\rho^{-4}).
	\end{equation}
	In particular, $\xi$ satisfies the bound
	\begin{equation}\label{apeq:difference-Jacobi}
		|J\xi|
		\le C\rho^{-2-2\epsilon_*}|\xi|
		+C\rho^{-\alpha}\|p(\rho,\cdot)\|_2+Ce^{-c\rho}.
	\end{equation}
	The implicit constants depend only on those in \eqref{apeq:coarse-hypothesis}, the fixed exponents, and the choices of nested tubes and cutoffs.
\end{lemma}

\begin{proof}
	Subtracting the two Allen--Cahn equations gives
\[
	(-\Delta+u^2+uv+v^2-1)w=0.
\]
Conjugation and localization yield \eqref{apeq:difference-equation}. The bound for $Q_\Lambda$ follows from the estimates $|u-H|+|v-H|\le C\rho^{-\alpha}e^{-\mu_0|t|}$ and $V_0=O(\rho^{-2})$ in the tube. The cutoff error $\mathcal C_\chi[\Lambda w]=[\mathcal L_0,\vartheta](\Lambda w)$ is supported where $\kappa\rho\le|t|\le\kappa_1\rho$. On this region, the bounds for $w$ and $\Lambda$, together with their first derivatives, give $|\mathcal C_\chi[\Lambda w]|\le Ce^{-c\rho}$ for some $c>0$.

We next compute $\sigma^{-1}\int_{\mathbb R}Q_\Lambda h^2\,dt$, the coefficient of $\xi$ obtained by projecting the potential term $Q_\Lambda\widetilde w$ onto $h$. Apply Lemma~\ref{lem:coarse} with larger cutoffs that equal one on a tube containing $\operatorname{supp}\widehat\vartheta$. On that tube, setting $u_1=u$ and $u_2=v$, we have
\[
	u_i=H+\xi_i h+(n-2)\rho^{-2}\eta+p_i,
	\quad \xi_i=O(\rho^{-\alpha}),\quad
	p_i=O_\mu(\rho^{-2\alpha}).
\]
Due to \eqref{apeq:geometric-expansion} and the exponential decay of $h$, the geometric term contributes
\begin{equation}\label{apeq:geometric-average}
	\frac1\sigma\int_{\mathbb R}\widehat\vartheta V_0h^2\,dt
	=-\frac{n-2}{2\rho^2}+O(\rho^{-4}).
\end{equation}
Writing $q_i=u_i-H$, we have
\[
	u^2+uv+v^2-3H^2
	=3H(q_1+q_2)+q_1^2+q_1q_2+q_2^2.
\]

The projection of $3(\xi_1+\xi_2)Hh$ against $\widehat\vartheta h^2$ vanishes by parity. The remaining leading term is $6(n-2)\rho^{-2}H\eta$, so \eqref{apeq:normal-identities} provides
\begin{equation}\label{apeq:profile-potential-average}
	\frac1\sigma\int_{\mathbb R}
	\widehat\vartheta(u^2+uv+v^2-3H^2)h^2\,dt
	=-\frac{n-2}{2\rho^2}+O(\rho^{-2\alpha}).
\end{equation}
Since $\widehat\vartheta=1$ on $|t|\le\kappa\rho$ and both $h^2$ and $H\eta h^2$ decay exponentially in $|t|$, their integrals against $1-\widehat\vartheta$ are $O(e^{-c\rho})$. Thus
\[
	\int_{\mathbb R}\widehat\vartheta h^2\,dt=\sigma+O(e^{-c\rho}),
	\qquad
	\int_{\mathbb R}\widehat\vartheta H\eta h^2\,dt=-\frac{\sigma}{12}+O(e^{-c\rho}).
\]
Adding \eqref{apeq:geometric-average} and \eqref{apeq:profile-potential-average} establishes \eqref{apeq:averaged-potential}.
	
Finally, we project \eqref{apeq:difference-equation} onto $h$ and use \eqref{apeq:projection-identities} to find
\[
J\xi
=\left(\frac{n-2}{\rho^2}
+\frac1\sigma\int_{\mathbb R}Q_\Lambda h^2\,dt\right)\xi
+\frac1\sigma\int_{\mathbb R}Q_\Lambda ph\,dt
-\frac1\sigma\int_{\mathbb R}\mathcal C_\chi[\Lambda w]h\,dt.
\]
The first coefficient is $O(\rho^{-2-2\epsilon_*})$ by \eqref{apeq:epsilon}. The Cauchy--Schwarz inequality and \eqref{apeq:difference-equation} bound the second term by $C\rho^{-\alpha}\|p\|_2$, and the last term is exponentially small. This proves \eqref{apeq:difference-Jacobi}.
\end{proof}

We now explain how to apply the orthogonal estimate at each step of Lemma~\ref{lem:fast-difference}. Suppose that $w=O_\mu(\rho^{-k})$ on the tube $\mathcal T_{\kappa,R}$, and let $\mu_0$ denote the exponential decay rate from the initial bounds on the larger tube $\mathcal T_{\kappa_3,R}$. Fix $0<\nu<\min\{\mu,\mu_0,\sqrt2\}$. Outside the inner tube, the initial bound implies the same polynomial decay with the smaller exponential rate $\nu$, since
\[
	\rho^{-\alpha}e^{-\mu_0|t|}
	\le C_k\rho^{-k}e^{-\nu|t|}
	\qquad\text{for }|t|\ge\kappa\rho.
\]
Consequently, $w=O_\nu(\rho^{-k})$ throughout $\mathcal T_{\kappa_3,R}$. Applying interior Schauder estimates to the equation $(-\Delta+u^2+uv+v^2-1)w=0$ on unit balls within this larger tube yields weighted local $C^{2,\gamma}$ bounds on $\mathcal T_{\kappa_2,R_1}$, provided $R$ is sufficiently enlarged. Multiplying by the cutoff $\vartheta\Lambda$ and extending by zero therefore gives
\[
	\|\widetilde w\|_{C^{2,\gamma}_{k,\nu}(D_{R_1})}<\infty.
\]

The initial bounds additionally imply that
\[
	\|Q_\Lambda\|_{C^{0,\gamma}(B_1(X_0)\cap D_{R_1})}
	\le C\rho_0^{-\alpha},
	\qquad \text{for } X_0=(\rho_0,t_0)\in D_{R_1}.
\]
Combining this estimate with the weighted bound for $\widetilde w$ yields $Q_\Lambda\widetilde w=O_{\nu,\gamma}(\rho^{-k-\alpha})$. The commutator of the cutoff function satisfies the same bound: on its support, we have $|t|\ge\kappa\rho$, and the residual exponential factor $e^{-(\mu_0-\nu)|t|}$ absorbs any fixed polynomial growth in $\rho$.

The projection operator $\Pi F=F-\sigma^{-1}h\int_{\mathbb R}Fh\,dt$ preserves these weighted $C^{0,\gamma}$ estimates, provided $\nu<\sqrt2$. Indeed, the integral term exhibits the same polynomial decay and local H\"older regularity in $\rho$, while $h$ itself decays exponentially at a rate of $\sqrt2$ in $|t|$. Consequently, we obtain
\[
	\mathcal L_0p
	=-\Pi(Q_\Lambda\widetilde w)+\Pi\mathcal C_\chi[\Lambda w]
	=O_{\nu,\gamma}(\rho^{-k-\alpha}).
\]
Furthermore, since $p=\Pi\widetilde w=O_\nu(\rho^{-k})$, Lemma~\ref{lem:orthogonal} implies that for any $0<\nu'<\nu$ and a sufficiently large choice of $R_2$,
\[
	\|p\|_{C^{2,\gamma}_{k+\alpha,\nu'}(D_{R_2})}<\infty.
\]
These resulting estimates are uniform whenever the initial bounds and the assumed $k$-decay bound are uniform. For solutions lying between fixed bounding minimizers, this uniformity propagates through the finitely many steps of Lemma~\ref{lem:fast-difference}.

\subsection{Differentiated equations and odd--even decomposition}
\label{acapp:derivatives-subsection}

Here, we collect the identities used in the monotonicity arguments of Section~\ref{sec:signs}, as well as the profile equation for the odd part required in Corollary~\ref{cor:sharp-profile}. Let us define
\[
	U^*(r,s)=U(s,r),\quad S=\frac{U-U^*}{2},\quad \text{and} \quad
	T=\frac{U+U^*}{2},
\]
and introduce the differential operator $K=s\partial_r-r\partial_s=-t\partial_\rho+\rho\partial_t$.
Finally, we set $\mathcal A=-\Delta+3(S^2+T^2)-1$.

\begin{lemma}\label{lem:derivatives}
Let $U$ be an invariant solution to the equation $-\Delta U=U-U^3$. For $r,s>0$, we have
\begin{equation}\label{apeq:radial-linearizations}
	\left(\mathcal L_U+\frac{m-1}{r^2}\right)U_r=0,
	\qquad
	\left(\mathcal L_U+\frac{m-1}{s^2}\right)U_s=0.
\end{equation}
The relevant commutators are given by
\begin{equation}\label{apeq:commutators}
	\begin{gathered}
		{}[\Delta,K]=(m-1)(r^{-2}+s^{-2})K,\\
		[\partial_\rho,\Delta]=A_\rho\partial_\rho+B_\rho\partial_t,
	\end{gathered}	
\end{equation}
where $A_\rho=-\frac{(n-2)(\rho^2+t^2)}{(\rho^2-t^2)^2}$ and
$B_\rho=\frac{2(n-2)\rho t}{(\rho^2-t^2)^2}$.
The odd and even parts, $S$ and $T$, respectively satisfy the equations
\begin{equation}\label{apeq:ST-equations}
	-\Delta S=S(1-S^2-3T^2),\quad
	-\Delta T=T(1-T^2-3S^2).
\end{equation}
Furthermore, we obtain the relations
\begin{equation}\label{apeq:angular-T-equation}
	\bigl(\mathcal A+(m-1)(r^{-2}+s^{-2})\bigr)(-KT)=6ST\,KS
\end{equation}
and
\begin{equation}\label{apeq:coupled-radial-equations}
	\begin{aligned}
		\left(\mathcal A+\frac{m-1}{rs}\right)S_\rho+6ST\,T_\rho
		&=\frac{(m-1)(r-s)}{\sqrt2\,r^2s^2}KS,\\
		\left(\mathcal A+\frac{m-1}{rs}\right)(-T_\rho)-6ST\,S_\rho
		&=-\frac{(m-1)(r-s)}{\sqrt2\,r^2s^2}KT.
	\end{aligned}
\end{equation}
Finally, the difference $q=S-H$ satisfies
\begin{equation}\label{apeq:odd-profile-equation}
	\mathcal L_0q
	=Bh+\left(A-\frac{n-2}{\rho}\right)q_\rho+Bq_t
	-3Hq^2-q^3-3ST^2.
\end{equation}
\end{lemma}

\begin{proof}
Differentiating the equation for $U$ yields \eqref{apeq:radial-linearizations}; see \cite[(4.5)--(4.6)]{Cabre2012}. This calculation relies only on the equation itself and the invariance of $U$. The identities in \eqref{apeq:commutators} follow directly from \eqref{apeq:laplacian}.
	
	Because both $S+T=U$ and $S-T=-U^*$ solve the Allen--Cahn equation, adding
	and subtracting their respective equations gives \eqref{apeq:ST-equations}.
	For a vector field $Y$ on the $(r,s)$-plane, we have
	\begin{equation}\label{apeq:differentiated-system}
		\mathcal A(YS)+6ST\,YT=[Y,\Delta]S,
		\quad
		\mathcal A(YT)+6ST\,YS=[Y,\Delta]T.
	\end{equation}
	Taking $Y=\mathcal K$ yields \eqref{apeq:angular-T-equation}. For
	$Y=\partial_\rho$, it is convenient to rewrite the commutator as
	\[
	[\Delta,\partial_\rho]f
	=\frac{m-1}{\sqrt2}\left(\frac{f_r}{r^2}+\frac{f_s}{s^2}\right)
	=\frac{m-1}{rs}f_\rho
	-\frac{(m-1)(r-s)}{\sqrt2\,r^2s^2}\mathcal Kf.
	\]
	Substituting this expression into \eqref{apeq:differentiated-system} directly leads to 
	\eqref{apeq:coupled-radial-equations}. Finally, subtracting the equation
	for $H$ from the first equation in \eqref{apeq:ST-equations} and applying
	\eqref{apeq:model-difference} yields \eqref{apeq:odd-profile-equation}.
\end{proof}

\bibliographystyle{abbrv}
\bibliography{allen_cahn_references}

\end{document}